%% file: massealf.tex
\documentclass[a4paper,11pt]{amsart}
\usepackage{amsmath}
\usepackage{amsmath, amsthm, amssymb, amsfonts}
\usepackage{amscd}
\usepackage[all]{xy}
\usepackage{graphicx,color}

\newtheoremstyle{Style}%
{.5em}{.5em}%
{\it}%
{}%
{\sc}%
{\ {\bf---}}%
{ }%
{}%

\newtheoremstyle{StyleRemarque}%
{\topsep}{.5em}%
{\it}%
{}%
{\slshape}%
{.}%
{ }%
{}%

\theoremstyle{Style}
\newtheorem{defn}{Definition}[section]
\newtheorem{lem}[defn]{Lemma}
\newtheorem{prop}[defn]{Proposition}
\newtheorem{thm}[defn]{Theorem}
\newtheorem{cor}[defn]{Corollary}
\newtheorem*{thm*}{Theorem}

\theoremstyle{StyleRemarque}

\newtheorem*{rem*}{Remark}

\newtheorem*{ex*}{Example}

\newcommand{\norm}[1]{\left\Vert#1\right\Vert}
\newcommand{\abs}[1]{\left\vert#1\right\vert}

\newcommand{\set}[1]{\left\{#1\right\}}

\newcommand{\appli}[3]{#1 \, : \, #2 \To #3}
\newcommand{\operator}[5]{ 
\begin{eqnarray*}
#1 :  \quad #2 &\To& #3 \\ 
#4 &\mapsto& #5 
\end{eqnarray*}}

\def\XXint#1#2#3{{\setbox0=\hbox{$#1{#2#3}{\int}$}
     \vcenter{\hbox{$#2#3$}}\kern-.5\wd0}}

\newcommand{\Rl}{\mathbb R}
\newcommand{\Cx}{\mathbb C}
\newcommand{\Ir}{\mathbb Z}
\newcommand{\Nl}{\mathbb N}
\newcommand{\Sph}{\mathbb S}
\newcommand{\Ball}{\mathbb B}

\newcommand{\To}{\longrightarrow}
\newcommand{\OO}{\mathcal O}

\newcommand{\DD}{\mathcal D}
\newcommand{\RR}{\mathcal R}
\newcommand{\XX}{\mathcal X}
\newcommand{\BB}{\mathcal B}
\newcommand{\ZZ}{\mathcal Z}

\DeclareMathOperator{\End}{End}
\DeclareMathOperator{\Ker}{Ker}
\DeclareMathOperator{\Ran}{Ran}
 
\DeclareMathOperator{\Tr}{Tr}
\DeclareMathOperator{\Ric}{Ric}
\DeclareMathOperator{\Scal}{Scal}

\DeclareMathOperator{\vol}{vol}

\DeclareMathOperator{\id}{id}

\DeclareMathOperator{\dirac}{\not{\hskip -1mm D}} 
\DeclareMathOperator{\Div}{div}

\begin{document}


\title{\bfseries{A mass for ALF manifolds.}}
\author{Vincent Minerbe}
\date{\today}

\maketitle

\begin{abstract}
We prove positive mass theorems on ALF manifolds, i.e. complete noncompact 
manifolds that are asymptotic to a circle fibration over a Euclidean base,
with fibers of asymptotically constant length.  
\end{abstract}


\section*{Introduction.}
In general relativity \cite{ADM}, one can define the mass of an asymptotically Euclidean spatial slice $(M^n,g)$ in a spacetime by the formula
\begin{equation}\label{AEmass}
\mu_g = \frac{1}{\omega_{n}}
\lim_{R \To \infty} \int_{\Sph^{n-1}(R)}  [  
\partial_j g_{i j} \, - \partial_j g_{i i}]   \iota_{\partial_i} dx.
\end{equation}  
Here, $\Sph^{n-1}(R)$ is the standard sphere with radius $R$ in $\Rl^n$, $\omega_n$ is the
volume of $\Sph^{n-1}(1)$. The terminology ``asymptotically Euclidean'' means $M$ minus a compact
subset is diffeomorphic to $\Rl^n$ minus a ball and the metric $g$ is asymptotic to the standard
metric on $\Rl^n$. The positive mass conjecture asserts this mass, if defined, is nonnegative when the scalar curvature $\Scal_g$ is nonnegative, vanishing only when $(M^n,g)$
is isometric to $\Rl^n$. It is a theorem in dimension $3 \leq n \leq 8$ (\cite{SY1,SY2}) and on spin manifolds of any dimension (\cite{Wit}). Moreover, R. Bartnik \cite{Bart} proved the mass is a genuine Riemannian invariant (under sharp assumptions). J. Lohkamp \cite{Loh} recently announced a proof of the positive mass theorem in all dimensions, without spin assumption. The interested reader should certainly have a look at \cite{LP,H1}. The mathematical interest of such a theorem is the rigidity result it involves: under a nonnegative curvature assumption, it asserts a model metric at infinity (the Euclidean metric, here) cannot be approached at any rate, the obstruction being precisely the mass.      

Positive mass theorems were proved in other settings. Asymptotically (real) hyperbolic manifolds were studied from this point of view in \cite{MO,AD,CH}, while asymptotically complex hyperbolic manifolds are treated in \cite{H2,BH}. Motivated by theoretical physics considerations, Xianzhe Dai \cite{Dai} pointed out a positive mass theorem for spin manifolds that are asymptotic to the product of a Euclidean space with a compact simply connected Calabi-Yau manifold. 

In this paper, we are interested in manifolds that are asymptotic to a circle fibration over a Euclidean base, with fibers of asymptotically constant length. For instance, there are complete Ricci flat metrics on $\Rl^2 \times \Sph^{n-2}$ that are asymptotic to the standard metric on $\Rl^{n-1} \times \Sph^1$. In dimension four, an example is provided by the Euclidean Schwarzschild metric:
$$
g = \frac{dr^2}{1 - \frac{2m}{r} } + r^2 d\omega^2 + \left( 1 - \frac{2m}{r} \right) dt^2.
$$ 
Another kind of example, adapted to the Hopf fibration at infinity, is the hyperk\"ahler Taub-NUT metric on $\Rl^4$: 
$$
g = \left( 1 + \frac{2m}{r} \right) \left( dr^2 + r^2 d\omega^2 \right) + \frac{1}{1 + \frac{2m}{r}} \eta^2,
$$
where $\eta$ is the standard contact form on the three-spheres (lots of details will be given about these formulas in the text). In these expressions, the parameter $m$ is bound to be nonnegative (so as to get a complete metric) and it is interpreted as a mass by physicists. It is tempting to ask for a positive mass theorem involving this $m$ in these examples. Indeed, some motivation comes from the main result of \cite{M2}: gravitational instantons with cubic volume growth are always asymptotic to circle fibrations over a Euclidean base (up to a finite group action) with an approximation rate of order $r^{-1}$, where $r$ is the distance to some point. Moreover, a class of Ricci flat metrics, including Schwarzschild examples enjoy similar properties (\cite{M2}). With this in mind, it seems interesting to wonder if there are non trivial examples of Ricci flat metrics $g$ that are closer to their model at infinity, namely for instance:
$$
g = g_{\Rl^3 \times \Sph^1} + \OO(r^{-1-\epsilon}),  \leqno{(A)}
$$
with $\epsilon>0$? The corresponding question with a nontrivial circle fibration at 
infinity can also be addressed. This note will prove in particular that the answer is no, because of a positive mass theorem. 

Let $\XX^{m+1}$ be the total space of a principal $\Sph^1$-bundle $\pi$ over 
$\BB^m:=\Rl^m \backslash \mathbb{B}^m$ (the exterior of the unit ball in $\Rl^m$, $m \geq 3$).
Given a positive number $L$, we introduce the vector field $T$ that is equal to  $\frac{L}{2\pi}$ times the infinitesimal generator of the $\Sph^1$ action and consider 
a ``model'' metric $h$, given by 
$$
h = \pi^*g_{\Rl^m} + \eta^2,
$$
where $\eta$ is a connection one-form, namely a $\Sph^1$ invariant one-form on $\XX$ such that $\eta(T)=1$. For such a one-form $\eta$, one may write $d\eta = \pi^* \omega$ for some two-form $\omega$ on the base $\BB$ and we will assume this ``curvature'' two-form $\omega$ decays at infinity: $\omega = \OO(r^{1-m})$ and $D\omega = \OO(r^{-m})$. Observe the fibers have length $L$. Metrics asymptotic to such metrics $h$ are known as ``ALF metrics''. 

Our main results are two positive mass theorems corresponding to such asymptotic models. 
The first positive mass theorem is valid for any circle fibration under a nonnegativity assumption on the Ricci curvature. The second one holds for asymptotic trivial fibrations, with a spin assumption and a nonnegativity assumption on the scalar curvature. We will denote by $B_R$ the preimage by $\pi$ of the ball of radius $R$ on the base. 

\begin{thm}
Let $(M^{m+1},g)$, $m \geq 3$, be a complete oriented manifold with nonnegative Ricci curvature. We assume  that, for some compact subset $K$, $M \backslash K$ is the total space of a circle fibration over $\Rl^m \backslash \mathbb{B}^m$, which can be endowed with a model metric $h$ such that   
$$
g = h + \OO(r^{2-m}), \quad \nabla^h g = \OO(r^{1-m}) \quad
\text{and} \quad \nabla^{h,2} g = \OO(r^{-m}).
$$
Then the quantity defined by
$$
\mu_{g,h}^{GB} = -\frac{1}{\omega_{m} L} \limsup_{R \To \infty} \int_{\partial B_R} *_{h} \; \Big{(} \Div_h g + d \Tr_h g -\frac{1}{2} d \; g(T,T) \Big{)}.
$$
is nonnegative and vanishes exactly when $(M,g)$ is the standard $\Rl^m\times \Sph^1$. Moreover, provided $\Ric_g$ is integrable, it is a Riemannian invariant (it depends only on $g$).
\end{thm}
Indeed, this quantity is the trace of a nonnegative quadratic form arising naturally in view of Witten's ideas. We stress the fact that the assumption $\Ric \geq 0$ is too strong for the asymptotically Euclidean setting (for Bishop theorem immediately implies the manifold is isometric to $\Rl^n$), while it allows many interesting examples when the volume growth is slower.  

\begin{thm}
Let $(M^{m+1},g)$, $m \geq 3$, be a complete \emph{spin} manifold with nonnegative scalar curvature. We assume there is a compact set $K$, a ball $B$ in $\Rl^{m}$ and a spin preserving diffeomorphism between $M \backslash K$ and $\Rl^{m} \backslash B \times \Sph^1$ such that  
$$
g = g_{\Rl^m \times \Sph^1} + \OO(r^{2-m}), \quad D g = \OO(r^{1-m}) \quad
\text{and} \quad D^{2} g = \OO(r^{-m}).
$$
Then the quantity defined by
$$
\mu_{g,h}^{D} = -\frac{1}{\omega_{m} L} \limsup_{R \To \infty} \int_{\partial B_R} *_{h_0} \; \Big{(} \Div_{h_0} g + d \Tr_{h_0} g \Big{)}.
$$
is nonnegative and vanishes exactly when $(M,g)$ is the standard $\Rl^m\times \Sph^1$. Moreover, provided $\Scal_g$ is integrable, it is a Riemannian invariant.
\end{thm}
 
There is a definite link with Dai's work \cite{Dai}. In his paper, he points out he cannot include Schwarzschild-like metrics in the discussion because the spin structure on the circles at infinity is not trivial. Indeed, physicists even built examples of complete spin manifolds with asymptotic (A), with nonnegative scalar curvature and with a \emph{negative} mass (\cite{BrH,CJ})! Here, we can ``justify'' the positivity of the mass of genuine Schwarzschild metrics, relying on the fact that their Ricci tensor vanishes (see section 4). Moreover, we can cope with nontrivial fibrations at infinity, which is interesting in view of the Taub-NUT example. 

The paper is organized as follows. In a first section, we introduce the class of metrics we are interested in. In a second section, we describe the analytical tools required for our arguments. We could have relied on Mazzeo-Melrose machinery; yet, we have chosen to include elementary proofs for everything we need. An advantage is we only make $C^2$ assumptions on the metric and never use complete Taylor expansions at infinity; the theory is thus simpler, closer to what is usually done in the asymptotically Euclidean case and it does not require any familiarity with pseudo-differential calculus. We emphasize nonetheless that the paper \cite{HHM} greatly inspired us, as well the text of some lectures given by Frank Pacard \cite{Pac}. In a third section, we prove the positive mass theorems. In the last section, we compute the masses in Schwarzschild and Taub-NUT examples; we also discuss Reissner-Nordstr\"om metrics, which have a negative mass.
 
We will often use an Einstein summation convention without mentionning it: when an index is repeated, the expression should be summed over this index (with a range depending on the context). 

\textbf{Acknowledgements.} I would like to thank Gilles Carron and Marc Herzlich for their interest in this work and for many interesting comments. This work benefited from the French ANR grant GeomEinstein.

\section{Metrics adapted to a circle fibration at infinity.}

Let us work on the total space $\XX^{m+1}$ of a principal circle fibration $\pi$ over 
$\BB:=\Rl^m \backslash \mathbb{B}^m$, $m \geq 3$. From a topological point of view, there are two families of such fibrations. In any dimension, we can consider the trivial circle fibration $\pi_0$ over $\BB=\Rl^m \backslash \Ball^m$. The total space is $\XX = \XX_0 := \Rl^m \backslash \Ball^m \times \Sph^1$. In dimension four (i.e. $m=3$), the Hopf fibration $\Sph^3 \To \Sph^2$ induces non trivial fibrations $\pi_k$, $k \in \Nl^*$, where the total space $\XX_k$ is the quotient of $\Rl^4 \backslash (2\Ball^4)$ by $\Ir_k$; the action of $\Ir_k$ is the complex scalar action of the $k$th unit root group on  $\Cx^2=\Rl^4$. 

We wish to define a class of ``model'' metrics $h$ on $\XX$. Their first feature is that the fibers of $\pi$ will have \emph{constant}
length $L$ ($0<L<+\infty$). It is therefore natural to introduce the vertical vector field $T$ that is equal to $\frac{L}{2\pi}$ times the infinitesimal generator of the $\Sph^1$-action; the flow of $T$ goes around the fibers in time $L$. As a second feature, $h$ should pullback 
the Euclidean metric on $\Rl^m$ as much as possible, i.e. on the orthogonal of the fibers. The way to do this is to pick a ``connection'' one-form 
$\eta$, namely a $\Sph^1$-invariant one-form $\eta$ such that $\eta(T)=1$. We define the corresponding model metric on $\XX$ by 
$$
h := \pi^*g_{\Rl^m} + \eta^2 = dx^2 + \eta^2,
$$ 
where $x_i := \pi^* \check{x}_i$, $1 \leq i \leq m$, denote the lifts of the canonical coordinates $\check{x}_i$ on the base $\BB \subset \Rl^m$.
In addition, we will require some decay for the curvature $d\eta$ of this connection. To express it, let us define $r := \sqrt{x_1^2 + \dots + x_m^2}$ and observe $d\eta = \pi^* \omega$ for some two-form $\omega$ on the base $\BB$ ($d\eta$ is $\Sph^1$-invariant and Cartan's magic formula $L_T \eta = \iota_T d\eta + d(\iota_T\eta)$ implies $\iota_T d\eta = 0$). We will assume
$$
\omega = \OO(r^{1-m}) \quad \text{and} \quad D\omega = \OO(r^{-m}).
$$
These estimates are easily seen to hold in the following examples.  

\vskip 0.5cm
 
\begin{ex*}[Trivial fibration.]
In the trivial fibration case, we can choose a trivialization $(x_1,\dots,x_m,e^{it})$ and put $\eta = \eta_0 :=dt$. Hence $d\eta_0=0$. Observe $\eta_0$ defines \emph{foliations} over $r^{-1}(R)$ for any $R$, corresponding to the product foliation of $\Rl^m \times \Sph^1$ by circles. The model metric is the flat $h_0:=dx^2+dt^2$.
\end{ex*}
\begin{ex*}[Hopf fibration.]
In the Hopf fibration case, we define connection forms $\eta_1$ on $\Rl^4=\Cx^2$ by the formula $\eta_k (z):= \abs{z}^{-2} g_{\Rl^4}(T,.)$, where $T$ is the infinitesimal generator of the scalar action of $\Sph^1 \subset \Cx^*$. Then $d\eta_k$ is the pull back of the volume form on $\Sph^2$ by $\pi_k$. Note $\eta_1$ is nothing but the standard 
contact form on $\Sph^3$ and $h_k := dx^2+\eta_k^2$ is the model at infinity for the multi-Taub-NUT metrics.
\end{ex*}

Let us describe a few features of such metrics. The coframe $(dx_1,\dots,dx_m,\eta)$ is
obviously orthonormal. We let $(X_1,\dots,X_m,T)$ be the dual frame. Uniqueness of the Levi-Civita connection ensures that, for horizontal vector fields $X$, $Y$  (namely $\eta(X)=\eta(Y)=0$),
$\pi_* \nabla_X^h Y = D_{\pi_*X} \pi_*Y$. The computation of brackets of horizontal vector fields, together with $T$, contains the geometric information about the fibration. First, since $\eta ([X_i,T]) = -d\eta(X_i,T) = 0$ and $\pi_* [X_i,T] = [\pi_*X_i, \pi_* T] = 0$, $T$ commutes with 
any vector field $X_i$. Besides, with $\eta ([X_i,X_j]) = -d\eta(X_i,X_j)$ and $\pi_* [X_i,X_j] = [\pi_*X_i, \pi_* X_j] = 0$, we obtain $[X_i,X_j] = -d\eta(X_i,X_j) T$. Koszul formula
$$
2 (\nabla_A^h B, C) = A \cdot (B,C) + B \cdot (A,C) - C \cdot (A,B)
                  - (A, [B,C])    - (B, [A,C])    + (C, [A,B])
$$
enables us to compute the Levi-Civita connection $\nabla^h$:
\begin{equation}\label{connect}
\nabla^h \eta = \frac{d\eta}{2} \qquad \text{and} \qquad
\nabla^h dx_i = \frac{\iota_{X_i} d\eta \otimes \eta + \eta \otimes \iota_{X_i} d\eta}{2}.
\end{equation}
In other words, the only non trivial Christoffel coefficients are
$$
h(\nabla_{X_i}^h T, X_j) = h(\nabla_{T}^h X_i, X_j) = - h(\nabla_{X_i}^h X_j, T) = \frac{d\eta(X_i,X_j)}{2}.
$$
In particular, we deduce:
\begin{equation}\label{estimconnexh}
\nabla^h X_i = \OO(r^{1-m}) \quad \text{and} \quad \nabla^h T = \OO(r^{1-m}).
\end{equation}
Moreover, since $\nabla^h_{X_i}X_i = \nabla^h_T T = 0$, we have for any function $u$ on $\XX$:
$$
\Delta_h u = -\Tr_h \nabla^h du 
= - (\nabla^h_{X_i} du)(X_i) - (\nabla^h_{T} du)(T) 
= - X_i \cdot X_i \cdot u - T \cdot T \cdot u.
$$
Given a local section of the circle bundle, one may consider a local vertical coordinate $t$ and work in the coordinates $(x_1,\dots,x_m,t)$. Observe $\partial_t = T$ and $\eta = dt + A_i dx_i$ for some functions $A_i$ independent of $t$. So we can
write $X_k = \partial_k - A_k \partial_t$ and therefore: 
\begin{equation}\label{exprlap}
\Delta_h = - \partial_{k k} - \partial_{t t} + 2 A_k \partial_{k t} - A_k^2 \partial_{t t}
+ \partial_k A_k \partial_t.
\end{equation}

\section{Analysis on asymptotic circle fibrations.}

In this section, we work on a complete Riemannian manifold $(M^{m+1},g)$, $m \geq 3$, such that for some compact subset $K$, $M \backslash K$ is diffeomorphic to the total space $\XX$ of a principal circle fibration $\pi$ over $\BB:=\Rl^m \backslash \mathbb{B}^m$, with:
$$
g = h + o(1) \quad \text{and} \quad \nabla^h g = o(r^{-1}),
$$
where $h$ is a model metric as in section 1. We intend to study the equation $\Delta_g u = f$ in weighted $L^2$ spaces.
The presentation is deeply influenced by \cite{HHM} and \cite{Pac}. 

Given a real number $\delta$ and a subset $\Omega$ of $M$, we define the weighted Lebesgue space
$$
L^2_\delta (\Omega) := \set{u \in L^2_{loc} \quad \Big{/} \quad 
\int_{\Omega \backslash K} u^2 r^{-2\delta} dvol_h < \infty}
$$
and endow it with the norm 
$$
\norm{u}_{L^2_\delta(\Omega)} := \left( \int_{\Omega \cap K} u^2 dvol_g 
+ \int_{\Omega \backslash K} u^2 r^{-2\delta} dvol_h \right)^{\frac{1}{2}}.
$$
Changing $K$ only produces equivalent norms. We will often write $L^2_\delta$ for 
$L^2_\delta (M)$. The following should be kept in mind: 
\begin{equation}\label{poidscroissance}
r^{a} \in L^2_{\delta}(M \backslash K) \Leftrightarrow \delta > \frac{m}{2} + a.
\end{equation} 
Note the Riemannian measures $dvol_g$ and $dvol_h$ on $M \backslash K$ can always be 
interchanged, thanks to our asymptotic assumption. For the same reason, the Riemannian 
connections $\nabla^g$ and $\nabla^h$ will be completely equivalent outside $K$. 

Any function $u \in L^2_\text{loc}(M \backslash K)$ can be written
$
u = \Pi_0 u + \Pi_\bot u
$
where $\Pi_0 u$ is obtained by computing the mean value of $u$ along the fibers: we set
$$
\Pi_0 u (x) = \frac{1}{L} \int_{\pi^{-1}(x)} u \, \eta. 
$$
It corresponds to a Fourier series decomposition along the fibers: $\Pi_0 u$ is the part 
in the kernel of $-T^2=-\partial_{t t}$; $\Pi_\bot u$ is the part in the positive eigenspaces of $-\partial_{t t}$.
The point is $\Delta_h$ commutes with the projectors $\Pi_0$ and $\Pi_\bot$: given a 
function $u$, we can use (\ref{exprlap}) to write locally on the base:
$$
\Pi_0 (\Delta_h u) (x) 
= \frac{1}{L} \int \left(- \partial_{k k} - \partial_{t t} + 2 A_k(x) \partial_{k t} - A_k(x)^2 \partial_{t t}
+ \partial_k A_k(x) \partial_t \right) u(x,t) dt,
$$
which simplifies into 
$$
\Pi_0 (\Delta_h u) (x) 
= -\frac{1}{L} \int \partial_{x x} u(x,t) dt 
= -\partial_{x x}  \frac{1}{L} \int u(x,t) dt 
= \Delta_h \left( \frac{1}{L} \int u(x,t) dt \right).
$$
This ensures $\Pi_0 \Delta_h = \Delta_h \Pi_0$ and also $\Pi_\bot \Delta_h = \Delta_h \Pi_\bot$.

Since this decomposition will prove useful, we introduce the following 
Hilbert spaces, depending on two real parameters $\delta$ and $\epsilon$:
$$
L^2_{\delta,\epsilon} (\Omega):= \set{u \in L^2_{loc}(\Omega) \quad \Big{/} \quad 
\norm{\Pi_0 u}_{L^2_\delta(\Omega \backslash K)} < \infty \quad \text{and} \quad 
\norm{\Pi_\bot u}_{L^2_\epsilon(\Omega \backslash K)} < \infty}.
$$
They are endowed  with the following Hilbert norm: 
$$
\norm{u}_{L^2_{\delta,\epsilon}(\Omega)} := \left( \norm{u}_{L^2(K\cap \Omega)}^2 
+ \norm{\Pi_0 u}_{L^2_\delta(\Omega \backslash K)}^2 
+ \norm{\Pi_\bot u}_{L^2_\epsilon(\Omega \backslash K)}^2 \right)^{\frac{1}{2}}.
$$
We also introduce the Sobolev space 
\begin{eqnarray*}
H^2_{\delta} := \Big{\{} u \in H^2_{loc}  \Big{/}
\norm{\nabla^g d\Pi_0 u}_{L^2_{\delta-2}(K^c)} + \norm{d \Pi_0 u}_{L^2_{\delta-1}(K^c)} + \norm{\Pi_0 u}_{L^2_\delta(K^c)} &<& \infty \\ 
\text{and} \quad \norm{\nabla^g d\Pi_\bot u}_{L^2_{\delta-2}(K^c)} + \norm{d \Pi_\bot u}_{L^2_{\delta-2}(K^c)} + \norm{\Pi_\bot u}_{L^2_{\delta-2}(K^c)} &<& \infty \Big{\}},
\end{eqnarray*}
endowed with the obvious Hilbert norm.   

In what follows, we will always write $A_R$ for the ``annulus'' defined by $R\leq r \leq 2R$ and $A^\tau_R$ for $2^{-\tau} R \leq r \leq 2^{\tau +1} R$ ($\tau \geq 0$). Similarly, the ``balls'' $K \cup \set{r \leq R}$ will be denoted by $B_R$. The letter $c$ will always denote a positive constant whose value changes from line to line. Its precise dependence on parameters will often be clear in the context; otherwise, we will write $c(\dots)$ to clarify it.

\subsection{A priori estimates.}

\subsubsection{A priori estimates on the kernel of $\Pi_\bot$.}

We aim at some a priori estimates for $\Delta_h$ on $\XX$. Let us begin with the kernel of $\Pi_\bot$: basically, we work with functions defined on exterior domain in $\Rl^m$ and the Laplace operator is the standard one! What follows is therefore standard (it is in \cite{Pac} and can be compared with \cite{LP,Bart,HHM}). We will nonetheless provide a few details of the proofs, in order to use them later.    

\begin{lem}\label{gardingbase}
Given $\delta \in \Rl$, there is a positive number $c=c(m,\delta)$ such that for any 
$R_0 \geq 1$ and any $u$ in $L^2_\delta \cap \Ker \Pi_\bot$, 
$$
\norm{\nabla^h du}_{L^2_{\delta-2}(B_{2R_0}^c)} + \norm{du}_{L^2_{\delta-1}(B_{2R_0}^c)} 
\leq c \left( \norm{\Delta_h u}_{L^2_{\delta-2}(B_{R_0}^c)} + \norm{u}_{L^2_\delta(B_{R_0}^c)} \right).
$$
\end{lem}

\proof
Scaling the usual Garding inequality for the Laplacian on $\Rl^{m}$ yields
$$
\norm{D du}_{L^2_{\delta-2} (A_R)}^2 + \norm{du}_{L^2_{\delta-1} (A_R)}^2 
\leq c \left( \norm{\Delta_h u}_{L^2_{\delta-2} (A^1_R)}^2 + \norm{u}_{L^2_\delta (A^1_R)}^2 \right).
$$
The formulas for the connection of $h$ imply 
$
\abs{\nabla^h du} \leq c \abs{D du} + c r^{1-m} \abs{du}
$
so we can replace $D$ by $\nabla^h$ in the estimate above. Using it for $R = 2^{i+1} R_0$, $i\in \Nl$, and summing over $i$, we get the desired inequality.
\endproof

To carry on, we need to use a $L^2$ spectral decomposition for 
the Laplace operator $\Delta_S$ on the unit sphere $\Sph^{m-1}$
in $\Rl^m$. Its eigenvalues are $\lambda_j := j (m-2+j)$, with $j \in \Nl$, and 
we denote by $E_j$ the corresponding eigenspaces. We also set  
$$
\delta_j := \frac{m}{2} + j
\quad \text{and} \quad
\nu_j^\pm := \frac{2-m}{2} \pm (\delta_j-1)
$$
for $j\in \Nl$. This simply means 
$\nu_j^+ = j$ and $\nu_j^- = 2 - m - j$. These numbers  $\nu_j^\pm$ are usually called 
``indicial roots''. Their basic property is the following: given an element $\phi_j$ 
of $E_j$, we have $\Delta \left( r^{\nu_j} \phi_j(\omega) \right) = 0$
outside $0$. It can easily be seen from the formula
\begin{equation}\label{lappol}
\Delta = -\partial_{rr} - \frac{m-1}{r} \partial_r + \frac{1}{r^2} \Delta_S -\partial_{t t}.
\end{equation}

\begin{defn}
We will say $\delta$ is critical if $\delta = \delta_j$ or $\delta = 2 - \delta_j$ for some 
$j$ in $\Nl$. 
\end{defn}

\begin{figure}[htb!]
\begin{center}
\input{expcrit.pstex_t}
\end{center}
\caption{The critical exponents.}
\end{figure}
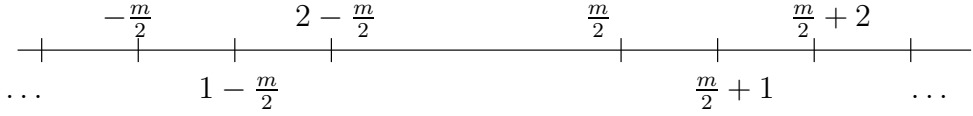

\begin{lem}\label{rellichbase}
If $\delta$ is not critical, there is a positive number $c=c(m,\delta)$ such that for any 
${R_0} \geq 1$ and any  $u$ in $L^2_\delta \cap \Ker \Pi_\bot$,
$$
\norm{u}_{L^2_{\delta}(B_{2{R_0}}^c)} 
\leq c \left( \norm{\Delta_h u}_{L^2_{\delta-2}(B_{{R_0}}^c)} 
+ \norm{u}_{L^2_{\delta}(A_{{R_0}})} \right).
$$
\end{lem}

\proof
If we first perform a spectral decomposition for the standard Laplace operator $\Delta_S$ along the sphere $\Sph^{m-1}$, our problem reduces to find an estimate on solutions $u_j \in L^2(\Rl_+, r^{m-1-2\delta}dr)$ to the radial equation 
$$
-\partial_{rr} u_j - \frac{m-1}{r} \partial_r u_j + \frac{\lambda_j}{r^2} u_j = f_j.
$$ 
Setting $r:=e^{s}$, $v_j(s):=u(r)$ and $g_j(s) =r^2 f_j(r)$, we turn the equation into
\begin{equation}\label{eqexp}
\frac{d}{ds} \left( v_j'(s)  e^{(m-2)s} \right) = (\lambda_j v_j(s) - g_j(s) ) e^{(m-2)s} 
\end{equation}
And we work in $L^2(\Rl_+, d\mu_\delta)$, with $d\mu_\delta :=e^{(m-2\delta)s}ds$. Given a truncature function $\chi$ vanishing on $B_{R_0}$ and equal to $1$ outside $B_{2 R_0}$, integration by parts provides
\begin{eqnarray*}
& & \int ((\chi v_j)')^2 d\mu_\delta \\
&=& \int \chi'^2 v_j^2 d\mu_\delta + \int (\chi^2 v_j)' \left(v'_j e^{(m-2)s} \right) e^{(2-2\delta)s} ds \\
&=& \int \chi'^2 v_j^2 d\mu_\delta - \int (\chi^2 v_j) (\lambda_j v_j - g_j) d\mu_\delta 
- (1-\delta) \int \chi^2 (v_j^2)' e^{(m-2\delta)s} ds \\
&=& \int \left[ \chi'^2  + (1-\delta) (\chi^2)' \right] v_j^2 d\mu_\delta 
+ \left[ (1-\delta)(m-2\delta) - \lambda_j \right]  \int (\chi v_j)^2 d\mu_\delta 
+ \int \chi^2 v_j g_j d\mu_\delta. 
\end{eqnarray*}
One might then use the Hardy inequality
$$
\frac{(m-2\delta)^2}{4} \int (\chi v_j)^2 d\mu_\delta \leq \int ((\chi v_j)')^2 d\mu_\delta 
$$
to deduce
$$
(\delta_j - \delta)(\delta_j + \delta -2) \int (\chi v_j)^2 d\mu_\delta
\leq \int \left[ \chi'^2  + (1-\delta) (\chi^2)' \right] v_j^2 d\mu_\delta 
+ \int \chi^2 v_j g_j d\mu_\delta.
$$
For all indices $j$ such that $(\delta_j - \delta)(\delta_j + \delta -2) >0$, 
we can use Young inequality do deduce
$$
\int_{2R_0}^\infty (\chi v_j)^2 d\mu_\delta
\leq c \int_{R_0}^{2R_0}  v_j^2 d\mu_\delta 
+ c \int \chi^2 g_j^2 d\mu_\delta, 
$$
which is what we need. For the finitely many indices $j$ such $2-\delta <\delta_j < \delta$, we consider
the function $w_j := \frac{1}{2 (1-\delta_j)}  \left( w_{j,+} - w_{j,-} \right)$, with
\begin{equation}\label{explicitform}
w_{j,\pm}(s) :=  \int_{R_0}^s e^{\nu_j^\pm(s-\sigma)} g_j(\sigma) d\sigma.
\end{equation}
Since $w_{j,\pm}' = \nu_j^\pm w_{j,\pm} + g_j$, integration by parts leads to
\begin{eqnarray*}
& &(m-2\delta) \int_{R_0}^{R_\infty} w_{j,\pm}^2(s) e^{(m-2\delta)s} \\
&=& -2 \nu_j^\pm \int_{R_0}^{R_\infty} w_{j,\pm}^2 d\mu_\delta - 2 \int_{R_0}^{R_\infty}  w_{j,\pm} g_j d\mu_\delta
+ \underbrace{w_{j,\pm}^2(R_\infty) e^{(m-2\delta)R_\infty}}_{\geq 0}.
\end{eqnarray*} 
Since $m-2\delta +2\nu_j^\pm$, it implies (Cauchy-Schwarz inequality):
$$
\int_{R_0}^\infty w_{j,\pm}^2 d\mu_\delta \leq c \int_{R_0}^\infty g_{j}^2 d\mu_\delta.
$$ 
Now, $w_j$ is a solution of (\ref{eqexp}) as well as $v_j$, so from EDO theory,
$w_j - v_j$ is a linear combination of $e^{\nu_j^\pm s}$. Hence its $L^2(d\mu_\delta)$-norm 
over $[R_0,+\infty[$ can be estimated by its $L^2(d\mu_\delta)$-norm over $[R_0,2R_0[$. It follows that $v_j$ can be estimated by
$$
\int_{R_0}^\infty v_{j}^2 d\mu_\delta \leq c \int_{R_0}^\infty g_{j}^2 d\mu_\delta
+ c \int_{R_0}^{2R_0} v_{j}^2 d\mu_\delta.
$$ 
The remaining case, $\delta <\delta_j < 2-\delta$ can be dealt with in a similar way. One defines
an explicit solution $w_j$ as above, replacing only $R_0$ by $+ \infty$ in formula (\ref{explicitform}). The integration by parts argument still works (because $w_{j,\pm}=o\left( e^{(\delta-\frac{m}{2})s}\right)$) and the fact that $v_j$ and $w_j$ belong to $L^2(\Rl_+,d\mu_\delta)$ forces them to coincide. 
\endproof

\vskip 0.5cm

\subsubsection{A priori estimates on the kernel of $\Pi_0$.}

The following lemma shows that estimates are basically \emph{better} on the kernel of $\Pi_0$. 
 
\begin{lem}\label{pasbase}
Given $\delta \in \Rl$ and a large number $R_0$, there is a constant $c$ such that for 
any function $u$ in $L^2_\delta \cap \Ker \Pi_0$, 
$$
\norm{\nabla^h du}_{L^2_{\delta}(B_{2R_0}^c)} 
+\norm{du}_{L^2_{\delta}(B_{2R_0}^c)} 
+ \norm{u}_{L^2_{\delta}(B_{2R_0}^c)} 
\leq c \left( \norm{\Delta_h u}_{L^2_{\delta}(B_{R_0}^c)} + \norm{u}_{L^2_{\delta}(A_{R_0})} \right).
$$
\end{lem}

\proof
Given parameters $R >>1$ and $0\leq \tau < \tau' \leq 1$, we can always choose a smooth nonnegative function $\chi$ that is equal to $1$ on $A_R^{\tau}$, vanishes outside $A_R^{\tau'}$, is $\Sph^1$-invariant and has gradient bounded by $c(\tau,\tau') R^{-1}$. The integration by part formula
$$
\int \abs{d(\chi u)}^2 dvol_h = \int \abs{d\chi}^2 u^2 dvol_h
+ \int \chi^2 u \Delta_h u \; dvol_h.
$$
can be used together with Young inequality to obtain
\begin{equation}\label{ippest1}
\norm{du}_{L^2(A_R^{\tau})}^2 \leq c(\epsilon) \norm{\Delta_h u}_{L^2(A_R^{\tau'})}^2 
+ \epsilon \norm{u}_{L^2(A_R^{\tau'})}^2,
\end{equation}
where $\epsilon$ can be chosen arbitrarily small provided $R$ is sufficiently large. 
Since $\Pi_0 u = 0$, we have
$$
\int_{\text{fiber}} (-T^2 u)u \; \eta = \int_{\text{fiber}} du(T)^2 \; \eta \geq c \int_{\text{fiber}} u^2\; \eta.
$$  
This implies $\norm{u}_{L^2(A_R^{\tau})} \leq c \norm{du}_{L^2(A_R^{\tau})}$
and therefore, with (\ref{ippest1}):
\begin{equation}\label{ippest2}
\norm{u}_{L^2(A_R^{\tau})}^2 + \norm{du}_{L^2(A_R^{\tau})}^2
\leq c(\epsilon) \norm{\Delta_h u}_{L^2(A_R^{\tau'})}^2 
+ \epsilon \norm{u}_{L^2(A_R^{\tau'})}^2.
\end{equation}
To get an order two estimate, we first write 
$$
\int \chi^2 \abs{\nabla^h du}^2 = \int (\nabla^h (\chi^2 du), \nabla^{h}du)
- 2 \int \chi (d\chi \otimes du,\nabla^h du).
$$
The Bochner Laplacian $\nabla^{h,*}\nabla^h$ and the Hodge Laplacian $\Delta_h = d d^*_h + d^*_h d$
only differ by the Ricci endomorphism, so an integration by parts provides: 
$$
\int \chi^2 \abs{\nabla^h du}^2 = \int \chi^2 (\Delta_h du,du) - \int \chi^2 \Ric_h(du,du)
- 2 \int \chi (d\chi \otimes du,\nabla^h du).
$$
Since the Hodge Laplacian commutes with $d$, another integration by part yields 
$$
\int \chi^2 (\Delta_h du,du) = \int (d \Delta_h u, \chi^2 du) = \int \chi^2 \abs{\Delta_h u}^2 - 2 \int \chi (d\chi,du) \Delta_h u.
$$
Putting these formulas alltogether and using Young inequality, we obtain:
$$
\int \chi^2 \abs{\nabla^h du}^2 \leq 
4 \int \chi^2 \abs{\Delta_h u}^2 + 6 \int \abs{d\chi}^2 \abs{du}^2 - 2 \int \chi^2 \Ric_h(du,du). 
$$
The upshot of this formula is the following estimate (observe $\Ric_h$ is bounded):
$$
\norm{\nabla^h du}_{L^2(A_R^\tau)}^2 \leq  c \norm{\Delta_h u}_{L^2(A_R^{\tau'})} 
+ c \norm{du}_{L^2(A_R^{\tau'})}^2. 
$$
With (\ref{ippest2}), we deduce:
$$
\norm{u}_{L^2(A_R^{\tau})}^2 + \norm{du}_{L^2(A_R^{\tau})}^2 + \norm{\nabla^h du}_{L^2(A_R^\tau)}^2
\leq c(\epsilon) \norm{\Delta_h u}_{L^2(A_R^{\tau'})}^2 + \epsilon \norm{u}_{L^2(A_R^{\tau'})}^2.
$$
We then set $\epsilon=0.5$, $\tau=0$, $\tau'=1$, we multiply these inequalities by $R^{-2\delta}$ and sum them for $R = 2^k R_0$, $k \in \Nl^*$ ($R_0$ is chosen large) to find: 
$$
\norm{u}_{L^2_\delta(B_{2R_0}^c)}^2 + \norm{du}_{L^2_\delta(B_{2R_0}^c)}^2 + \norm{\nabla^h du}_{L^2_\delta(B_{2R_0}^c)}^2 \leq  
c \norm{f}_{L^2_\delta(B_{R_0}^c)}^2 + c \norm{u}_{L^2_\delta(A_{R_0})}^2.
$$
\endproof

\begin{rem*}
When the fibration is trivial, one can also study the integral kernel of the resolvent,
as in \cite{Dav}. Using a Fourier decomposition, we are left to prove that for positive $k$, the integral kernel
$$
P(x,y) := \left( \frac{1+d(o,y)}{1+d(o,x)} \right)^{\delta} (\Delta_{\Rl^m} + k^2)^{-1} (x,y)
$$
defines a bounded operator on $L^2(\Rl^m)$. Indeed, it is true for instance on a complete manifold with $\Ric \geq 0$ and $\forall t\geq 1, \, A t^\nu \leq \vol B(x,t) \leq B t^\nu$ for some $\nu >2$. Noticing
$$
(\Delta_{\Rl^m} + k^2)^{-1} = \frac{1}{\sqrt{\pi}} \int_0^\infty \frac{ e^{-t\Delta} e^{-k^2 t}}{\sqrt{t}} dt
\quad \text{and} \quad 
\frac{1+d(o,y)}{1+d(o,x)} \leq 1 + d(x,y),
$$
we can use Li-Yau Gaussian estimate on the heat kernel \cite{LY} to get
$$
P(x,y) \leq c (1 + d(x,y)) \int_0^\infty \frac{e^{-\frac{d(x,y)^2}{c t}}}{\vol B(x,\sqrt{t})} \frac{e^{-k^2 t}}{\sqrt{t}} dt.
$$
It implies $P(x,y) \leq c e^{-\frac{d(x,y)}{c}}$ and therefore $\int P(x,y) dx \leq c$ and
$\int P(x,y) dy \leq c$, which is enough to ensure the boundedness of $P$ on $L^p$, $1 < p < \infty$. 
\end{rem*}

\vskip 0.5cm

\subsubsection{The main estimate.}

\begin{prop}\label{estimee}
If $\delta$ is not critical, there is a constant $c$ and a compact set $B$ such that 
for any $u$ in $L^2_{\delta,\delta-2}$,
$$
\norm{u}_{H^2_\delta} 
\leq c \left( \norm{\Delta_g u}_{L^2_{\delta-2}} + \norm{u}_{L^2(B)} \right).
$$
\end{prop}

\vskip 0.5cm

\proof
Since $\Delta_h$ commutes with $\Pi_0$ and $\Pi_\bot$, we can write the equation $\Delta_g u = f$ outside a large ball $B_{R_0}$ as follows:
$$
\left\lbrace
\begin{array}{l}
\Delta_h u_0 = f_0 + \Pi_0 (\Delta_h - \Delta_g) u \\
\Delta_h u_\bot = f_\bot + \Pi_\bot (\Delta_h - \Delta_g) u
\end{array}
\right.
$$
We denote $\Pi_0 u$ by $u_0$, $\Pi_\bot u$ by $u_\bot$, etc, to make the equations easier to read.
We apply lemmata \ref{gardingbase} and \ref{rellichbase} to the first equation and lemma \ref{pasbase}
to the second equation, which results in   
\begin{eqnarray*}
& & \norm{\nabla^h du_0}_{L^2_{\delta-2} (B_{2{R_0}}^c)} + \norm{du_0}_{L^2_{\delta-1} (B_{2{R_0}}^c)} 
+ \norm{u_0}_{L^2_\delta (B_{2{R_0}}^c)}\\
&+& \norm{\nabla^h du_\bot}_{L^2_{\delta-2} (B_{2{R_0}}^c)} + \norm{du_\bot}_{L^2_{\delta-2} (B_{2{R_0}}^c)} 
+ \norm{u_\bot}_{L^2_{\delta-2} (B_{2{R_0}}^c)}\\
&\leq& c \norm{f}_{L^2_{\delta-2} (B_{{R_0}}^c)} + c \norm{u}_{L^2_\delta (B_{2{R_0}})} 
+ c \norm{(\Delta_h - \Delta_g) u}_{L^2_{\delta-2}(B_{{R_0}}^c)}.
\end{eqnarray*}
Since
\begin{equation}\label{difflap}
\abs{(\Delta_h - \Delta_g)u} 
\leq c \abs{g -h} \abs{ \nabla^h du} + c \abs{\nabla^g - \nabla^h } \abs{du},
\end{equation}
there is a function $\epsilon$ going to zero at infinity such that
$$
\norm{(\Delta_h - \Delta_g) u}_{L^2_{\delta-2}(B_{{R_0}}^c)}
\leq \norm{\epsilon(r) \nabla^h du}_{L^2_{\delta-2} (B_{{R_0}}^c)}
+ \norm{r^{-1} \epsilon(r) du}_{L^2_{\delta-2} (B_{{R_0}}^c)}.
$$
Therefore, $\norm{(\Delta_h - \Delta_g) u}_{L^2_{\delta-2}(B_{{R_0}}^c)}$
can be bounded by
\begin{eqnarray*}
&\epsilon({R_0})& \Big{(} \norm{\nabla^h du_0}_{L^2_{\delta-2} (B_{2{R_0}}^c)} 
+  \norm{\nabla^h du_\bot}_{L^2_{\delta-2} (B_{2{R_0}}^c)}  
+  \norm{du_0}_{L^2_{\delta-1} (B_{2{R_0}}^c)} \Big{)} \\
+ &\epsilon({R_0})&  \Big{(} \norm{du_\bot}_{L^2_{\delta-2} (B_{2{R_0}}^c)}   
 +  \norm{\nabla^h du}_{L^2_{\delta-2} (A_{{R_0}})}
+ \norm{du}_{L^2_{\delta-2} (A_{{R_0}})} \Big{)}. 
\end{eqnarray*}
Using this and choosing ${R_0}$ large enough, we find
\begin{eqnarray*}
& & \norm{\nabla^h du_0}_{L^2_{\delta-2} (B_{2{R_0}}^c)} + \norm{du_0}_{L^2_{\delta-1} (B_{2{R_0}}^c)}
+ \norm{u_0}_{L^2_\delta (B_{2{R_0}}^c)}\\
&+& \norm{\nabla^h du_\bot}_{L^2_{\delta-2} (B_{2{R_0}}^c)} + \norm{du_\bot}_{L^2_{\delta-2} (B_{2{R_0}}^c)} 
+ \norm{u_\bot}_{L^2_{\delta-2} (B_{2{R_0}}^c)}\\
&\leq& c \norm{f}_{L^2_{\delta-2} (B_{{R_0}}^c)} + c \norm{u}_{L^2 (B_{2{R_0}})} 
+ c \norm{\nabla^h d u}_{L^2(B_{2{R_0}})} + c \norm{du}_{L^2(B_{2{R_0}})}.
\end{eqnarray*}
Owing to the asymptotic of the metric, $\nabla^h$ can be changed into $\nabla^g$ 
in this estimate. Since the standard Garding inequality provides
\begin{eqnarray*}
& & \norm{\nabla^g du}_{L^2 (B_{2{R_0}})} + \norm{du}_{L^2 (B_{2{R_0}})}
\leq c \norm{f}_{L^2 (B_{4{R_0}})} + c \norm{u}_{L^2 (B_{4{R_0}})}, 
\end{eqnarray*}
we can conclude:
$$
\norm{u}_{H^2_\delta}\\
\leq c \norm{f}_{L^2_{\delta-2,\delta-2}} + c \norm{u}_{L^2 (B_{4{R_0}})}.
$$
\endproof

\subsection{Mapping properties.}

\subsubsection{Fredholmness.}

We would like to investigate some properties of the unbounded operator 
\operator{P_\delta}{\mathcal{D}(P_\delta)}{L^2_{\delta-2,\delta-2}}{u}{\Delta_g u} 
whose domain $\mathcal{D}(P_\delta)$ is the dense subset of $L^2_{\delta,\delta-2}$ 
whose elements $u$ have their (distribution) Laplacian in $L^2_{\delta-2,\delta-2}$. 
Proposition \ref{estimee} has a direct consequence. In view of Rellich's theorem and 
Peetre's lemma (cf. \cite{LiM}, p. 171), it implies $\Ker P_\delta$ is finite dimensional
(for any $\delta$) and $\Ran P_\delta$ is closed (for any noncritical $\delta$).
The usual $L^2$ pairing identifies the topological dual space of $L^2_{\delta,\delta-2}$
(resp. $L^2_{\delta-2,\delta-2}$) with $L^2_{-\delta,2-\delta}$ (resp. $L^2_{2-\delta,2-\delta}$). For this identification, the adjoint $P^*_\delta$ of $P_\delta$ is nothing but
\operator{P^*_\delta}{\mathcal{D}(P^*_\delta)}{L^2_{-\delta,2-\delta}}{u}{\Delta_g u} 
where the domain $\mathcal{D}(P^*_\delta)$ is the dense subset of $L^2_{2-\delta,2-\delta}$ 
whose elements $u$ have their (distribution) Laplacian in $L^2_{-\delta,2-\delta}$. 
Observing $\Ker P^*_\delta \subset \Ker P_\eta$ for some large non critical $\eta$, we see that
$\Ker P^*_\delta$ is always finite dimensional. We have proved the  

\begin{prop}
If $\delta$ is not critical, then $P_\delta$ is Fredholm. 
\end{prop}

This discussion also provides a first surjectivity result.

\begin{prop}\label{surj1}
If we further assume $\Ric_g \geq - c r^{-2}$, then $P_\delta$ is surjective for any noncritical $\delta > 2 - \frac{m}{2}$.  
\end{prop}

\begin{rem*}
The assumption on $\Ric_g$ holds for instance when $\nabla^{h,2} g = \OO(r^{-2})$. Corollary \ref{surj2} below will provide another surjectivity result, without such an order two estimate on $g$, but with a slightly stronger assumption up to order one.  
\end{rem*}

\proof
For noncritical $\delta$, $P_\delta$ has closed range, so it is surjective as
soon as $P_\delta^*$ is injective. Since $\Ric_g \geq - c r^{-2}$, we can use \cite{SC} (or even \cite{M1} for a stronger statement) to ensure the following family of scaled Sobolev inequalities: with $n=m+1$, any smooth function $u$ with compact support in $A_R$ satisfies 
\begin{equation}\label{sobolev}
\norm{u}_{L^{\frac{2n}{n-2}}} \leq c \frac{R}{(\vol B_R)^\frac{1}{n}} \norm{du}_{L^2}.
\end{equation}  
With a Moser iteration (\cite{BKN} for instance), we deduce for any harmonic function $u$:
$$
\norm{u}_{L^{\infty}(\partial B_R)} \leq \frac{c}{(\vol B_R)^{\frac{1}{2}}} \norm{u}_{L^2(A^1_R)}.
$$  
If $\delta>2-\frac{m}{2}$, Cauchy-Schwarz inequality then implies that any harmonic 
function $u \in L^2_{2-\delta}$ goes to zero at infinity and therefore vanishes 
(maximum principle). So $P_\delta^*$ is injective as soon as $\delta>2-\frac{m}{2}$.
It follows that $P_\delta$ is surjective for any noncritical $\delta>2-\frac{m}{2}$.
\endproof

\subsubsection{Solving an exterior problem.}

We first solve the equation $\Delta_h u = f$ on an exterior domain.

\begin{lem}\label{resolutionmodele}
Given a noncritical $\delta$ and a large number $R_0$, we can define a bounded operator
$
\appli{G_h}{L^2_{\delta-2}(B_{R_0}^c)}{H^2_{\delta}(B_{R_0}^c)} 
$ 
such that $\Delta_h \circ G_h = \id$. 
\end{lem}

\proof
We will define $G_h$ in three steps, relying on an orthogonal decomposition of $L^2_{\delta-2}(B_{R_0}^c)$.  

For the first step, we pick a function $f$ in $L^2_{\delta-2} \cap \Ker \Pi_\bot$ and we assume it has no component 
along the eigenspaces $E_j$ of $\Delta_S$ such that $(\delta_j - \delta)(\delta_j + \delta -2) >0$. Given $R >> R_0$, we can use 
standard elliptic theory to solve the equation $\Delta_h u_{R} = f$ in $L^2(B_R \backslash B_{R_0})$, with Dirichlet boundary condition. Adapting the proof of lemma \ref{rellichbase} (on $B_R \backslash B_{R_0}$, with $\chi=1$), we obtain
\begin{equation}\label{estimext}
\norm{u_{R}}_{L^2_{\delta-2}} \leq c \norm{f}_{L^2_{\delta-2}},
\end{equation}
with $c$ independent of $R$. Elliptic regularity then bounds the $H^2$ norm of $u_R$ over compact subsets in terms of $\norm{f}_{L^2_{\delta-2}}$ , so that we can use Rellich theorem and a diagonal argument to extract a sequence $u_R$ converging to a function $u$ in $H^1_{loc}$, with $\Delta_h u = f$
and $u=0$ on $\partial B_{R_0}$. Taking a limit in (\ref{estimext}), one gets $\norm{u}_{L^2_{\delta-2}(B_{R_0}^c)} \leq c \norm{f}_{L^2_{\delta-2}}$.
From (\ref{gardingbase}) (plus standard elliptic arguments near $\partial B_{R_0}$, we deduce an estimate on the derivatives:
\begin{equation}\label{h2delta}
\norm{u}_{H^2_\delta(B_{R_0}^c)} \leq c \norm{f}_{L^2_{\delta-2}(B_{R_0}^c)}.
\end{equation}
We need to show that such a $u$ is uniquely defined,  i.e. independent of the choice 
of extracted sequence. The difference $v$ between two such functions $u$ obeys 
$\Delta_{\Rl^m} v =0$ so, from EDO theory, its modes read
$
v_j(r,\omega) = r^{\nu_j^+} \phi_j^+(\omega) + r^{\nu_j^-} \phi_j^-(\omega)
$
with $\phi_j^\pm$ in $E_j$. Since $v_j$ is in $L^2_\delta$ and vanishes on the boundary, 
we have $\phi_j^\pm=0$, so $v=0$. We can therefore set $G_h f := u$.

As a second step, we observe the same approach can be used for a function $f$ in $L^2_{\delta-2} \cap \Ker \Pi_0$. We can still obtain the functions $u_R$ by solving the same problem and the proof of lemma \ref{pasbase} can be adapted to provide the estimate (\ref{estimext}), provided $R_0$
is large enough. This makes it possible to extract a converging subsequence as above, 
yielding a function $u$ such that $\Delta_h u = f$, $u=0$ on $\partial B_{R_0}$ and
satisfying (\ref{h2delta}). As for unicity, we consider the difference $v$ 
between two such functions $u$: $v$ is in $L^2_{\delta,\delta-2}\cap \Ker \Pi_0$,
vanishes on $\partial B_{R_0}$ and obeys $\Delta_h v =0$. Given a large number $R$, we can 
choose a cutoff function $\chi$ that is equal to $1$ on $A_R$, vanishes on $A^1_R$ and 
has gradient bounded by $10/R$. Then
$$
\int \abs{d(\chi v)}^2 = \int \abs{d\chi}^2 v^2 + \int \chi^2 v \Delta_h v
$$
implies
$
\norm{dv}_{L^2_{\eta}(A_R)} \leq c \norm{v}_{L^2_{\eta-1}(A^1_R)}
$
for any exponent $\eta$. Since $\Pi_0 v=0$, we get  
$$
\norm{v}_{L^2_{\eta}(A_R)} \leq c \norm{dv}_{L^2_{\eta}(A_R)} 
\leq c \norm{v}_{L^2_{\eta-1}(A^1_R)} \quad \text{ for any exponent } 
\eta.
$$
This improves $L^2_\eta$ estimates into $L^2_{\eta-1}$ estimates. So $v \in L^2_\delta$ 
implies $v \in L^2_\eta$ for any $\eta$. In particular, $v$ is in $L^2$ and thus vanishes. 
So $u$ is well defined and we can set $G_h f := u$.

As a third step, we consider those $f$ in $L^2_{\delta-2} \cap \Ker \Pi_\bot$ whose only component in the spectral decomposition
of $\Delta_S$ is in $E_j$, with $(\delta_j - \delta)(\delta_j + \delta - 2) <0$. In case $2-\delta<\delta_j <\delta$, we set 
$G_h f := u$, with
\begin{equation}\label{formexpl}
u(r,\omega) := \frac{1}{2 (1-\delta_j)}  \left(r^{\nu_j^+} \int_{R_0}^r t^{1-\nu_j^+} f(t,\omega) dt 
- r^{\nu_j^-} \int_{R_0}^r t^{1-\nu_j^-} f(t,\omega) dt  \right),
\end{equation}
The proof of lemma \ref{rellichbase} yields
\begin{equation}\label{est1bas}
\norm{u}_{L^2_{\delta}(B_{R_0}^c)} \leq c  \norm{f}_{L^2_{\delta-2}(B_{R_0}^c)}.
\end{equation}
With lemma \ref{gardingbase}, we obtain:
\begin{equation}\label{est2bas}
\norm{u}_{H^2_{\delta}(B_{2R_0}^c)} 
\leq c  \norm{f}_{L^2_{\delta-2}(B_{R_0}^c)}.
\end{equation}
Since $u$ vanishes on $\partial B_{R_0}$, standard elliptic estimates 
improve (\ref{est2bas}) into
$$
\norm{u}_{H^2_{\delta}(B_{R_0}^c)} 
\leq c  \norm{f}_{L^2_{\delta-2}(B_{R_0}^c)}.
$$
The case $\delta <\delta_j < 2-\delta$ is dealt with similarly, replacing $R_0$ by $+\infty$ in (\ref{formexpl}). 
Note $u$ does not vanish on $\partial B_{R_0}$ in this setting. But it can be bounded by 
$c \norm{f}_{L^2_{\delta-2}(B_{R_0}^c)}$ on $\partial B_{R_0}$; since $\Delta_S u_j = \lambda_j u_j$, 
$\norm{u}_{H^2(\partial B_{R_0})}$
is then bounded by $c \norm{f}_{L^2_{\delta-2}(B_{R_0}^c)}$, so the argument above still works.
\endproof

A perturbation argument extends this result to a more general setting. 

\begin{prop}\label{resolution}
Given a noncritical $\delta$ and a large number $R_0$, we can define a bounded operator
$
\appli{G_g}{L^2_{\delta-2}(B_{R_0}^c)}{H^2_{\delta}(B_{R_0}^c)} 
$
such that $\Delta_g \circ G_g = \id$. 
\end{prop}

\proof
Thanks to lemma \ref{resolutionmodele}, we can write
$
\Delta_g  = \Delta_h \left[ \id + G_h \left( \Delta_g - \Delta_h \right)  \right].
$   
With (\ref{difflap}), we obtain: 
$
\norm{(\Delta_h - \Delta_g) u}_{L^2_{\delta,\delta-2}(B_{R_0}^c)} \leq \epsilon(R_0) \norm{u}_{H^2_\delta(B_{R_0}^c)}.
$
Since $G_h$ is bounded from $L^2_{\delta,\delta-2}(B_{R_0}^c)$ to $H^2_\delta(B_{R_0}^c)$, we deduce 
$G_h \left( \Delta_g - \Delta_h \right)$ defines a bounded operator on $H^2_\delta(B_{R_0}^c)$, with 
norm strictly inferior to $1$ for $R_0$ large enough. So
$\id + \left( \Delta_g - \Delta_h \right) \Delta_h^{-1}$ is an automorphism of $H^2_\delta(B_{R_0}^c)$ and 
$
G_g := \left[ \id + G_h \left( \Delta_g - \Delta_h \right)  \right]^{-1} G_h
$
is a bounded operator from $L^2_{\delta,\delta-2}(B_{R_0}^c)$ to $H^2_\delta(B_{R_0}^c)$, 
with $\Delta_g G_g = \Delta_h G_h = \id$.
\endproof 

This lemma can be used to build functions which are harmonic outside a compact set and have some prescribed asymptotics.

\begin{cor}\label{harmprescrites}
Let us further assume
$$
g = h + \OO(r^{-a}) \quad \text{and} \quad \nabla^h g = \OO(r^{-a-1}) \quad \text{for some }a>0.
$$
Then, given $j \in \Nl$ and $\phi \in E_j$, there are functions 
$\aleph_{j,\phi}^\pm$ that are harmonic outside a compact set and can be written $\aleph_{j,\phi}^\pm = r^{\nu_j^\pm} \phi + v_\pm$ with $v_+$ in $H^2_\eta$ for any $\eta > \delta_j-a$ and $v_-$ in $H^2_\eta$ for any $\eta > 2-\delta_j-a$.
\end{cor}

\proof
For noncritical $\delta > \delta_j$, since $r^{\nu_j^+} \phi$ is in $L^2_\delta$ and $\Delta_h (r^{\nu_j^+} \phi ) = 0$, lemma \ref{gardingbase} implies 
$r^{\nu_j^+} \phi$ is in $H^2_\delta$. With (\ref{difflap}), we deduce
$\Delta_g (r^{\nu_j^+} \phi) \in r^{-a} L^2_{\delta-2} = L^2_{\delta-a-2}$. 
Now we can use lemma \ref{resolution} to solve $\Delta_g u = - \Delta_g ( r^{\nu_j^+} \phi )$ outside a compact set and put $\aleph_{j,\phi}^+:= \chi ( r^{\nu_j^+} \phi + u )$ for some smooth nonnegative function $\chi$ which vanishes on a large compact set and is equal to $1$ outside a larger compact set. The construction of $\aleph_{j,\phi}^-$ follows the same lines.
\endproof

\subsubsection{Decay jumps.}

The following lemma is the key to understand the growth of solutions to our equations. 

\begin{lem}\label{sautsmodele}
Consider an equation $\Delta_h u = f$ with $u$ in $L^2_{\delta}(B_{R_0}^c)$ and $f$ in $L^2_{\delta'-2}(B_{R_0}^c)$ for noncritical exponents $\delta > \delta'$ and a large number $R_0$. Then there is an element $v$ of
$L^2_{\delta',\delta'-2}(B_{R_0}^c)$ such that $u-v$ is a linear combination of 
the following functions:
\begin{itemize}
\item $r^{\nu_j^+} \phi_j$ with $\phi_j$ in $E_j$ and $\delta' < \delta_j < \delta$; 
\item $r^{\nu_j^-} \phi_j$ with $\phi_j$ in $E_j$ and $\delta' < 2-\delta_j < \delta$.
\end{itemize}
\end{lem}

\proof
We will build $v$ step by step, starting from the solution $\tilde{v}$ of $\Delta_h \tilde{v} = f$ provided by lemma \ref{resolutionmodele}: $\tilde{v}$ is in $L^2_{\delta'-2,\delta'}(B_{R_0}^c)$. Consider $w:=\Pi_0 (u - \tilde{v})$ and look at its modes $w_j$. The equation $\Delta_h w_j = 0$ implies 
$
w_j = r^{\nu_j^+} \phi_j^+ + r^{\nu_j^-} \phi_j^-
$ 
with $\phi_j^\pm$ in $E_j$. Observing
$r^{\nu_j^+} \in L^2_\eta \Leftrightarrow \eta > \delta_j$ and 
$r^{\nu_j^-} \in L^2_\eta \Leftrightarrow \eta > 2-\delta_j$, one can see that each term is either  in 
$L^2_{\delta'-2,\delta'}$, so that we can add it to $\tilde{v}$ and forget it, or satisfies the conditions 
in the statement. What about $z:=\Pi_\bot (u - \tilde{v})$ ? This $z$ satisfies 
$\Delta_h z = 0$, $\Pi_0 z=0$ and $z \in L^2_\delta$. As in the proof of lemma \ref{resolutionmodele}, 
for any $\eta$, $L^2_\eta$ estimates on $z$ improve into $L^2_{\eta-1}$ estimates, so $z \in L^2_\delta$ 
implies $z \in L^2_\eta$ for any $\eta$ and the proof is complete. 
\endproof

Let us generalize this.

\begin{prop}\label{sauts}
We again assume
$$
g = h + \OO(r^{-a}) \quad \text{and} \quad \nabla^h g = \OO(r^{-a-1}) \quad \text{for some }a>0.
$$
Let us consider an equation $\Delta_g u = f$ with $u$ in $L^2_{\delta}(K^c)$ and $f$ in $L^2_{\delta'-2}(K^c)$ for noncritical exponents $\delta > \delta'$. Then, up to enlarging 
$K$, there is an element $v$ of $L^2_{\delta',\delta'-2}(K^c)$ such that $u-v$ is a linear combination of the following functions:
\begin{itemize}
\item $\aleph_{j,\phi_j}^+$ with $\phi_j$ in $E_j$, if $\delta' < \delta_j < \delta$; 
\item $\aleph_{j,\phi_j}^-$ with $\phi_j$ in $E_j$, if $\delta' < 2-\delta_j < \delta$.
\end{itemize} 
\end{prop}

\proof
Thanks to lemma \ref{estimee}, the equation $\Delta_g u = f$, with $u$ in $L^2_{\delta,\delta-2}$ and $f$ in $L^2_{\delta'-2}$, ensures $u \in H^2_\delta$. With (\ref{difflap}), we obtain
$$
\Delta_h u \in L^2_{\delta'-2} + L^2_{\delta-a-2}.
$$ 
So if we pick any noncritical $\eta \geq \max(\delta',\delta-a)$, we have $\Delta_h u \in L^2_{\eta-2}$. Lemma \ref{sautsmodele} then says $u$ admits a decomposition
$$
u = u_1 + \sum_{j_+} r^{\nu_{j_+}^+} \phi_{j_+}^+ + \sum_{j_-} r^{\nu_{j_-}^-} \phi_{j-}^-
$$
where $u_1$ belongs to $L^2_{\eta,\eta-2}$, $\eta < \delta_{j_+} < \delta$ and 
$2- \delta < 2-\delta_{j_-} < 2-\eta$. With lemma \ref{harmprescrites}, we can therefore write 
$$
u = u_2 + \sum_{j_+} \aleph_{j_+, \phi_{j_+}^+}^+ + \sum_{j_-} \aleph_{j_-, \phi_{j_-}^-}^-
$$
with $u_2$ in $L^2_{\eta,\eta-2}$,  $\eta < \delta_{j_+} < \delta$ and 
$2- \delta < 2-\delta_{j_-} < 2-\eta$. 

If $\delta-a \leq \delta'$, we are done. 

If not, observe $\Delta_g u_2 = f$ outside a compact set and $u_2$ belongs to $L^2_{\eta,\eta-2}$. 
So we can repeat the argument with $u_2$ in the role of $u$ and $\eta$ in the role of $\delta$. 
In a finite number of steps, we are in the first case.
\endproof

We can deduce another surjectivity result (compare with \ref{surj1}).

\begin{cor}\label{surj2}
If we assume
$$
g = h + \OO(r^{-a}) \quad \text{and} \quad \nabla^h g = \OO(r^{-a-1}) \quad \text{for some }a>0,
$$
then $P_\delta$ is surjective for any noncritical value 
$\delta > 2- \frac{m}{2}$.  
\end{cor}

\proof
For any noncritical $\delta \geq 2$, $\Ker P_\delta^*$ consists of harmonic $L^2$ functions hence is trivial, so $P_\delta$ is surjective . Now, choose some 
$\delta'$ in $]2- \frac{m}{2}, 2[$ and pick some function $f$ in $L^2_{\delta'-2}$.
In particular, $f$ is in $L^2_{2-2}$, so there is a solution $u \in L^{2}_{2,2-2}$ 
of $\Delta_g u = f$. The proposition \ref{sauts}, with $\delta = 2$ and $\delta' = \delta'$, implies $u \in L^{2}_{\delta',\delta'-2}$, for there is no critical exponent in $]\delta',2[$ ! $P_{\delta'}$ is therefore surjective.  
\endproof

\subsection{An extension to the Dirac operator.}

If we endow $\Rl^m \times \Sph^1$ with its trivial spin structure, we can define
a Dirac operator $\dirac_{h_0}$ whose square $\dirac_{h_0}^2$ acts diagonally
as the Laplace operator $\Delta_{h_0}$ in a constant trivialization. As a 
consequence, the a priori estimates that we have proved for $\Delta_{h_0}$ 
are also available for $\dirac_{h_0}^2$. We wish to work on a complete spin 
Riemannian manifold $(M^{m+1},g)$, $m \geq 3$, such that 
for some compact subset $K$, $M \backslash K$ is spin-diffeomorphic to 
$\Rl^m \backslash \mathbb{B}^m \times \Sph^1$, with:
$$
g = h_0 + o(1), \quad \nabla^{h_0} g = o(r^{-1})
\quad \text{and} \quad \nabla^{h_0,2} g = o(r^{-2}).
$$
A perturbation argument easily yields the 

\begin{prop}\label{estimeedirac}
If $\delta$ is not critical, there is a constant $c$ and a compact set $B$ such that 
for any $\psi$ in $L^2_{\delta,\delta-2}$,
$$
\norm{\psi}_{H^2_\delta} 
\leq c \left( \norm{\dirac_g^2 \psi}_{L^2_{\delta-2}} + \norm{\psi}_{L^2(B)} \right).
$$
\end{prop}

The functional spaces are defined in the obvious way, using a constant trivialization. The 
proof is basically the same as that of proposition \ref{estimee}. Note however the estimate on
the second derivative of $g$: we need this to control the difference between the model Dirac
Laplacian $\dirac_{h_0}^2$ and the operator $\dirac_g^2$, and more precisely the $0$th order 
term (which of course vanishes for the Laplace operator on functions):
$$
\abs{\dirac_g^2 \psi - \dirac_{h_0}^2 \psi} \leq c\abs{g -h_0} \abs{ \nabla^{h_0,2} \psi} + c \abs{\nabla^{h_0} g} \abs{\nabla^{h_0} \psi} + c \left( \abs{\nabla^g - \nabla^{h_0} }^2
+ \abs{\nabla^{h_0,2} g} \right) \abs{\psi}.
$$

The Fredholmness of the corresponding operator $P_\delta$, $\delta$ noncritical, follows immediately. The decay jump phenomenons carry over to this setting in exactly the same way. 
In particular, the proof of corollary \ref{surj2} yields the
\begin{cor}\label{surjdirac}
If $\Scal_g \geq 0$, $P_\delta$ is surjective for any noncritical value 
$\delta > 2- \frac{m}{2}$.  
\end{cor}

\proof
Lichnerowicz formula and the assumption $\Scal_g \geq 0$ ensure the $L^2$ kernel of $\dirac$ is trivial,
so we can adapt the corresponding proof for the Laplace operator.  
\endproof

\section{Towards a mass.}

In this section, we introduce a notion of mass for our metrics. The first paragraph contains a few 
(classic) algebraic computations which are useful in the sequel. The second and third 
paragraph develop two points of view corresponding to the two standard examples of Dirac type
operators. 

\subsection{Algebraic preliminaries.}

Let $M^n$ be a Riemannian manifold and let $E$ be a bundle of left modules over the Clifford algebra bundle $Cl(TM)$. For definitions and basic facts about spin geometry, we refer to \cite{LM}. We assume $E$ is endowed with a compatible Euclidean metric $(.,.)$ and metric connection $\nabla$, whose curvature tensor is $R$. This data determines a Dirac type operator $\DD$ on $E$ and a section $\RR$ of $\End E$. In a local orthonormal basis $(e_a)_a$, these are 
$$
\DD = \sum_{a=1}^n e_a \cdot \nabla_{e_a}
\quad \text{and} \quad
\RR = \frac{1}{4} \sum_{a, b=1}^n [e_a \cdot ,e_b \cdot ]  R_{e_a,e_b}. 
$$

\begin{rem*}
Commutators are easier to handle than Clifford products, for the latter are not antisymmetric. The identity
\begin{equation}\label{comprod}
[e_a \cdot ,e_b \cdot ]   = 2 ( \delta_{a b} + e_a \cdot e_b \cdot)   
\end{equation}
makes the translation. Moreover, brackets are skew-symmetric with respect to the Euclidean metric: \begin{equation}\label{comzero}
([e_a \cdot ,e_b \cdot ]   \psi, \psi) = 0.
\end{equation}
 \end{rem*}

Given sections $\alpha$ and $\beta$ of $E$, we define a one-form $\zeta_{\alpha,\beta}$ on $M$ by the following formula (as in \cite{AD,H2,Dai}): 
$$
\zeta_{\alpha,\beta} (X) := (\nabla_X \alpha + X \cdot \DD \alpha, \beta ).
$$
The point is:
$$
d^* \zeta_{\alpha,\beta} =  (\DD \alpha, \DD \beta) - (\nabla \alpha, \nabla \beta) 
- (\RR \alpha , \beta).
$$ 
Integrating this Lichnerowicz-type formula over a domain $\Omega$, we get
\begin{eqnarray*}
\int_\Omega \left[ (\nabla \alpha, \nabla \beta) + (\RR \alpha , \beta) 
-  (\DD \alpha, \DD \beta) \right] dvol 
= - \int_\Omega d^* \zeta_{\alpha,\beta} dvol 
= \int_\Omega d * \zeta_{\alpha,\beta}
\end{eqnarray*}
and Stokes formula provides
\begin{equation}\label{ipp}
\int_\Omega \left[ (\nabla \alpha, \nabla \beta) + (\RR \alpha , \beta) 
- \int_\Omega (\DD \alpha, \DD \beta) \right] dvol 
=  \int_{\partial \Omega} * \zeta_{\alpha,\beta}.
\end{equation}

As a technical device, we consider the two-form $\omega_{\alpha,\beta}$ defined by
$$
\omega_{\alpha,\beta}(X,Y) = ([X\cdot,Y\cdot] \alpha, \beta) 
$$
and observe $d^* \omega_{\alpha,\beta}= 4\zeta_{\alpha,\beta} - 4\zeta_{\beta,\alpha}$,
so that Stokes theorem implies:
\begin{equation}\label{permut}
 \int_{\partial \Omega} * \zeta_{\alpha,\beta}
=  \int_{\partial \Omega} * \zeta_{\beta,\alpha}.
\end{equation}

\vskip 0.5cm

\subsection{The Gauss-Bonnet case.}

Let us consider a complete oriented Riemannian manifold $(M^{m+1},g)$, $m \geq 3$, that fits into the setting of section 1: for some compact subset $K$, $M \backslash K$ is diffeomorphic to the total space of a circle fibration $\pi$ over $\Rl^m \backslash \mathbb{B}^m$, which we endow with a model metric $h = dx^2 + \eta^2$. We assume the metric $g$ is asymptotic to $h$ in the following sense:
$$
g = h + \OO(r^{2-m}), \quad \nabla^h g = \OO(r^{1-m}) \quad
\text{and} \quad \nabla^{h,2} g = \OO(r^{-m}).
$$

In this paragraph, we work on the exterior bundle $\Lambda M (=:E)$, endowed with the Levi-Civita connection $\nabla^g$ and we use the Gauss-Bonnet operator $d+d^*$ as Dirac-type operator $\DD$. The Clifford product $v \cdot := \epsilon_v - \iota_{v}$ is obtained from the exterior product $\epsilon_v := (v,.) \wedge$ and the interior product $\iota_v$. In this setting, $\RR$ preserves the form degree. In degree $1$, it reduces to the natural action of the Ricci tensor. 

Let us introduce the vector space $\ZZ$ spanned by $X_1,\dots,X_m$. As a first step toward the definition of a mass, we wish to build $g$-harmonic one-forms that are asymptotic to any element of the dual vector space $\ZZ^*$. We fix a small positive number $\epsilon$ (say $0 < \epsilon < 1/2$).

\begin{lem}\label{constructharm}
For any element $\widetilde{\alpha}$ in $\ZZ^*$, there is a one-form $\alpha$ on $M$ such that $(d+d^*) \alpha =0$ and $\alpha = \widetilde{\alpha} + \beta$, with $\beta = \OO(r^{2-m+\epsilon})$ and $\nabla^g \beta \in L^2_{1 - \frac{m}{2} + \epsilon}$.
\end{lem}

\proof It is enough to prove the claim for $\widetilde{\alpha}:=dx_k$. Let us choose a truncature function $\chi$ that vanishes on a large compact set and is equal to $1$ outside a larger compact set. Since $\Delta_g x_k$ is in $L^2_{\delta'-2}$ for any $\delta' > 3 - \frac{m}{2}$, we can apply weighted analysis (corollary \ref{surj1} or \ref{surj2}) to construct a function $u_k$ in $H^2_{\delta'}$ such that $\Delta_g u_k = - \Delta_g \left( \chi x_k \right)$, for any noncritical $\delta' > 3 - \frac{m}{2}$. Setting $\beta_k := du_k$, we therefore obtain $(d+d^*) \beta_k = - (d+d^*) d\left( \chi x_k \right)$, with $\beta_k \in L^2_{\delta}$ and $\nabla \beta_k \in L^2_{\delta-1}$, for any $\delta > 2 - \frac{m}{2}$ such that $\delta'=\delta +1$ is not critical. We may choose $\delta = 2 - \frac{m}{2} + \epsilon$. Then the one-form $\alpha_k := d\left( \chi x_k \right) + \beta_k$ obeys $(d+d^*) \alpha_k =0$. Since the scaled Sobolev inequality (\ref{sobolev}) holds and $\Delta_g \beta_k$ is equal to $\Delta_g dx_k =\OO(r^{-m})$ outside a compact set, a Moser iteration argument (\cite{GT}; more precisely, lemma 3.9 in \cite{TV} or lemma A.3 in \cite{M2}) can be coupled with the $L^2$ estimate to provide a $L^\infty$ estimate: $\beta_k = \OO(r^{2-m+\epsilon})$.
\endproof

Now, the basic formula we need is given by (\ref{ipp}) and (\ref{permut}): for any one-form $\alpha$ as in lemma \ref{constructharm} and any large number $R$, we have 
\begin{equation}\label{formule}
\int_{B_R} \Big{[} \abs{\nabla^g \alpha}^2 + \Ric_g(\alpha,\alpha) \Big{]} dvol 
=  \int_{\partial B_R} * \zeta_{\widetilde{\alpha},\widetilde{\alpha}} 
+ 2 \int_{\partial B_R} * \zeta_{\widetilde{\alpha},\beta} 
+ \int_{\partial B_R} * \zeta_{\beta,\beta}. 
\end{equation}
Our next aim consists in understanding the asymptotic behaviour of the right hand-side as much as possible. The mass is to be the limit of this quantity when the domains $B_R$ are larger and larger. Let us tackle the second and third terms on the right-hand side.   

\begin{lem}\label{smallterms}
There is a sequence $(R_i)_i$ going to infinity such that 
$$
\lim_{i \To \infty} \int_{\partial B_{R_i}} * \zeta_{\widetilde{\alpha},\beta} = 0
\quad \text{and} \quad
\lim_{i \To \infty} \int_{\partial B_{R_i}} * \zeta_{\beta,\beta} = 0
$$
\end{lem}

\proof Since $\beta$ is in $ L^2_{2 - \frac{m}{2} + \epsilon}$ and $\nabla \beta$ is in $ L^2_{1 - \frac{m}{2} + \epsilon}$, one can find a sequence $R_i$ going to infinity such that
$\displaystyle{
\left( \int_{\partial \Omega_{R_i}} \left( \abs{\beta}^2 + r^2 \abs{\nabla^g \beta}^2 \right) 
dvol \right)^{\frac{1}{2}} = o(R_i^{2 - \frac{m}{2} + \epsilon - \frac{1}{2}}).
}$
The lemma then follows from Cauchy-Schwarz inequality and $\nabla^g \widetilde{\alpha} = \OO(r^{1-m})$.
\endproof
Now, we need to compute the first term in the right-hand side of (\ref{formule}). Observe $h$ identifies $\ZZ$ and $\ZZ^*$, so any $Z$ in $\ZZ$ corresponds to a well defined $\widetilde{\alpha}= \widetilde{\alpha}_Z$ in $\ZZ^*$.
\begin{lem}\label{leadingterm} 
We have  
$$
\zeta_{\widetilde{\alpha},\widetilde{\alpha}} 
= - ( \Div_h g ) (Z) \widetilde{\alpha}
-\frac{1}{2} \left[ d (\Tr_h g) (Z) \,  \widetilde{\alpha} + d \left( g(Z,Z) \right) \right] 
+ \OO(r^{3-2m}).  
$$
\end{lem}

\proof 
Let us $g$-orthonormalize the frame field $(X_1,\dots,X_m,T)$ into $(e_a)_a$ and set 
$\omega_{c d} := g(\nabla^g e_c, e_d) - g(\nabla^h e_c, e_d)$. Since the connection on the \emph{tangent} bundle $TM$ reads $\nabla^g = \nabla^h + \omega_{c d} \; e^c \otimes e_d$, the connection on the \emph{cotangent} bundle satisfies:
$$
\nabla^g \widetilde{\alpha} = \nabla^h \widetilde{\alpha} - \omega_{d c} \; \widetilde{\alpha}(e_c) e^d.
$$
Therefore, we can write $\zeta_{\widetilde{\alpha},\widetilde{\alpha}} = \xi_1 + \xi_2$, with
\begin{eqnarray*}
\xi_1 &=& \frac{1}{2} g([e_a\cdot ,e_b \cdot]  \nabla_{e_b}^h \widetilde{\alpha} , \widetilde{\alpha} ) e^a, \\
\xi_2 &=& - \frac{1}{2} \omega_{d c}(e_b) \, \widetilde{\alpha}(e_c) g([e_a\cdot ,e_b \cdot]   e^d, \widetilde{\alpha} ) e^a.
\end{eqnarray*}
We are lead to estimate terms like
$$
 g([e_a\cdot ,e_b \cdot]  e^d, \widetilde{\alpha} )
=  g([\epsilon_{e_a} - \iota_{e_a},\epsilon_{e_b} - \iota_{e_b}] e^d, \widetilde{\alpha} ) \\
= - g([\iota_{e_a}, \epsilon_{e_b}] e^d + [\epsilon_{e_a},\iota_{e_b}] e^d, \widetilde{\alpha}) 
$$
(recall $\widetilde{\alpha}$ is a one-form). A little algebra provides
$
[\epsilon_{e_a},\iota_{e_b}] e^d = 2  \delta_{b d} e^a - \delta_{a b} e^d,
$
which leads to:
\begin{equation}\label{coefs}
\frac{1}{2} g([e_a\cdot ,e_b \cdot]  e^d, \widetilde{\alpha} )
=  \delta_{a d} \widetilde{\alpha}(e_b) - \delta_{b d} \widetilde{\alpha}(e_a).   
\end{equation}
We can use (\ref{coefs}), the asymptotic of $g$ and (\ref{connect}) to compute $\xi_1$:
\begin{eqnarray*}
\xi_1 &=& g\left( \nabla^h_{e_b} \widetilde{\alpha}, e^d \right) \left[ \delta_{a d} \widetilde{\alpha}(e_b) - \delta_{b d} \widetilde{\alpha}(e_a) \right] e^a\\
&=& \nabla^h_{Z} \widetilde{\alpha} - g\left( \nabla^h_{e_b} \widetilde{\alpha}, e^b \right) \widetilde{\alpha}\\
&=& \nabla^h_{Z} \widetilde{\alpha} + d^*_h \widetilde{\alpha} + \OO(r^{3-2m}) \\
&=& \OO(r^{3-2m}).
\end{eqnarray*}
As for $\xi_2$, we use (\ref{coefs}) again:
$
\xi_2 
=  \omega_{d c}(e_b) \, \widetilde{\alpha}(e_c) \left[ \delta_{b d} \widetilde{\alpha}(e_a) - \delta_{a d} \widetilde{\alpha}(e_b) \right] e^a. 
$
Since Koszul formula provides
\begin{eqnarray*}
2 \omega_{x y}(e_w) 
&=& - g(e_w,[e_x,e_y]) - g(e_x,[e_w,e_y]) + g(e_y,[e_w,e_x])- 2 g(\nabla^h_{e_w} e_x, e_y) \\
&=& - g(e_w,\nabla^h_{e_x} e_y) - g(e_x,\nabla^h_{e_w} e_y) + g(e_y,\nabla^h_{e_w} e_x)
+ g(e_w,\nabla^h_{e_y} e_x) \\
& &+ g(e_x,\nabla^h_{e_y} e_w) - g(e_y,\nabla^h_{e_x} e_w) - 2 g(\nabla^h_{e_w} e_x, e_y) \\
&=& (\nabla^h_{e_x}g) (e_y,e_w)  - (\nabla^h_{e_y}g) (e_x,e_w) + (\nabla^h_{e_w}g) (e_x,e_y),
\end{eqnarray*}
we eventually find:
\begin{eqnarray*}
\xi_2
&=& \frac{1}{2} \left[ 
(\nabla^h_{e_b}g) (Z,e_b) \,  \widetilde{\alpha} - (\nabla^h_{Z}g) (e_b,e_b) \,  \widetilde{\alpha} + (\nabla^h_{e_b}g) (e_b,Z) \,  \widetilde{\alpha} \right]\\
& & - \frac{1}{2} \left[  
(\nabla^h_{e_a}g) (Z,Z) \,  e^a - (\nabla^h_{Z}g) (e_a,Z) \,  e^a 
+ (\nabla^h_{Z}g) (e_a,Z) \,  e^a
\right] \\
&=& -\frac{1}{2} \left[ 
2 ( \Div_h g ) (Z) \widetilde{\alpha} + d (\Tr_h g) (Z) \,  \widetilde{\alpha} 
+ d \left( g(Z,Z) \right)  -2 h(\nabla^h Z,Z) 
\right]  + \OO(r^{3-2m}) \\
&=& -\frac{1}{2} \left[ 
2 ( \Div_h g ) (Z) \widetilde{\alpha} + d (\Tr_h g) (Z) \,  \widetilde{\alpha} 
+ d \left( g(Z,Z) \right)   
\right]  + \OO(r^{3-2m}).  
\end{eqnarray*}
\endproof

The computations are done, it is time to draw a theorem, which requires a definition.
\begin{defn}
On $M \backslash K$, we can define a one-form $q_{g,h}$ with values in the quadratic forms on $\ZZ$ by the formula
$$
q_{g,h}(Z) =  
- ( \Div_h g ) (Z) \widetilde{\alpha}_Z
-\frac{1}{2} \left[ d (\Tr_h g) (Z) \,  \widetilde{\alpha}_Z + d \left( g(Z,Z) \right) \right] 
.
$$
The ``mass'' quadratic form $\mathcal{Q}_{g,h}$ is the quadratic form defined on $\ZZ$ by:
$$
\mathcal{Q}_g(Z) := \frac{1}{\omega_{m} L}
\limsup_{R \To \infty} \int_{\partial B_R} *_{h} \; q_{g,h}(Z) 
$$
where $\omega_{m}$ is the volume of the standard $\Sph^{m-1}$ and $L$ is the asymptotic length of fibers.
\end{defn}

Why this normalization constant? The factor $1/L$ is there to make the mass independent of the length of the asymptotic circles. The normalization by the volume of a sphere is more anecdotic. 
This corresponding positive mass theorem is the following.

\begin{thm}\label{GBmass}
Let $(M^{m+1},g)$, $m \geq 3$, be a complete oriented manifold with nonnegative Ricci curvature. We assume  that, for some compact subset $K$, $M \backslash K$ is the total space of a circle fibration over $\Rl^m \backslash \mathbb{B}^m$, which can be endowed with a model metric $h$ such that   
$$
g = h + \OO(r^{2-m}), \quad \nabla^h g = \OO(r^{1-m}) \quad
\text{and} \quad \nabla^{h,2} g = \OO(r^{-m}).
$$
Then $\mathcal{Q}_{g,h}$ is a nonnegative quadratic form. It vanishes exactly when $(M,g)$ is the standard $\Rl^m\times \Sph^1$.
\end{thm}

\proof
Formula (\ref{formule}), together with $\Ric \geq 0$, lemma \ref{smallterms} and 
lemma \ref{leadingterm}, provides
\begin{equation}\label{massineq}
\int_{M} \Big{[} \abs{\nabla^g \alpha_Z}^2 + \Ric_g(\alpha_Z,\alpha_Z) \Big{]} dvol \leq \omega_{m} L \; \mathcal{Q}_{g,h}(Z).
\end{equation}
Since $\Ric \geq 0$, we deduce $\mathcal{Q}_{g,h}(Z) \geq 0$. 

Now we assume $\mathcal{Q}_{g,h}=0$. In view of (\ref{massineq}), the $g$-harmonic one-forms $\alpha_Z$ are then $g$-parallel. We have therefore built $m$ $g$-parallel one-forms $\alpha_1,\dots,\alpha_{m}$ that are asymptotic to $dx_1,\dots,dx_m$. We also put $\alpha_{m+1} := *_g (\alpha_1 \wedge \dots \wedge \alpha_{m})$, so as to obtain $m+1$ $g$-parallel one-forms that are linearly independent outside a compact set, hence linearly independent on the whole $M$. This yields a parallel coframe field on $M$: $(M,g)$ is flat. And the only flat manifold with the required asymptotic is $\Rl^m\times \Sph^1$.
\endproof

Note $\mathcal{Q}_{g,h}$ vanishes if and only if its trace vanishes. So 
$\mu_{g,h}^{GB} := \Tr \mathcal{Q}_{g,h}$ plays the role of a numerical mass. We have:
\begin{eqnarray}\label{formule1}
\mu_{g,h}^{GB} = -\frac{1}{\omega_{m} L}
\limsup_{R \To \infty} \int_{\partial B_R} *_{h} \; \Big{(} \Div_h g + d \Tr_h g -\frac{1}{2} g(T,T)\Big{)}.
\end{eqnarray}
In this formula, the integrand can be replaced by the Hodge star of
$$
\left[
(-\nabla^h_{X_j}g)(X_i,X_j)  + d\left( g(X_j,X_j) \right)(X_i) - (\nabla^h_{T}g)(X_i,T) + \frac{d\left( g(T,T) \right)(X_i)}{2} \right] dx_i.
$$
This can be simplified for, if we write 
$$
(\nabla^h_{T} g)(X_i,T)
= T \cdot g(X_i,T) - g(X_i,\nabla^h_{T}T)
- h(\nabla^h_{T}X_i,T) + (h-g)(\nabla^h_{T} X_i,T),
$$
we can use (\ref{connect}) and the closeness of $g$ to $h$ to get
\begin{eqnarray*}
(\nabla^h_{T} g)(X_i,T) *_h dx_i
&=& T \cdot g(X_i,T) *_h dx_i + \OO(r^{3-2m}) \\
&=& d \Big{(} g(X_i,T) *_h (dx_i \wedge \eta) \Big{)} + X_i \cdot g(X_i,T) *_h\eta + \OO(r^{3-2m}), 
\end{eqnarray*}
which integrates to zero at infinity! In the same spirit, we have
$$
(\nabla^h_{X_j}g)(X_i,X_j) = X_j \cdot g(X_i,X_j) + \OO(r^{3-2m}),
$$
so we are left with
\begin{eqnarray*}
\mu_{g,h}^{GB} = \frac{1}{\omega_{m} L}
\limsup_{R \To \infty} \int_{\partial B_R} *_{h} \; \Big{(}   \left[
 X_j \cdot g(X_i,X_j)  -  X_i \cdot g(X_j,X_j)  \right] dx_i 
 - \frac{d\left( g(T,T)\right)}{2}   \Big{)}.
\end{eqnarray*}

The term between brackets is similar to the usual expression of the mass in the asymptotically Euclidean setting (see (\ref{AEmass}) in the introduction). It is the contribution from the base.
More precisely, one can average the metric along the fibers and compute the mass of the asymptotically Euclidean metric induced on the base: it is this term. The other term is related to the variation of the length of fibers.  

We now turn to the geometric invariance of the mass: does the mass really depends on $h$ or is it a Riemannian invariant of $g$ ? 

\begin{prop}\label{welldef}
Provided the Ricci tensor is in $L^1$, the mass $\mu_{g,h}^{GB}$ depends only on $g$ and not on $h$: it is a Riemannian invariant. 
\end{prop}

\proof
Let us consider two metrics $h=dx^2+\eta^2$ and $h'=dx'^2 + \eta'^2$ built from two circle fibrations with connections as above. The proof of the unicity of the mass is in three steps (compare \cite{Bart,CH}). We first prove that, given an adapted (i.e. as above) $h$-orthonormal frame field $(X_1,\dots,X_m,T)$, we can find an adapted $h'$-orthonormal frame field $(X'_1,\dots,X'_m,T')$ that is $r^{2-m}$-close to $(X_1,\dots,X_m,T)$ and whose covariant derivatives are $r^{1-m}$-close. Secondly, we show that the computation of the mass term does not depend on the choice of ``spheres at infinity'' ($\partial B_R = r^{-1}(R)$ or $\partial B'_R = r'^{-1}(R)$). And finally, we prove that the difference between the integrands corresponding to $h$ and $h'$ is the sum of an exact form and of a negligible term. 

The first step will be completed once we have proved the infinitesimal generators $T$ and $T'$ of both circle actions are $r^{-1}$-close. To do so, we consider the smooth loop $\sigma$ defined on $M$ by $\sigma(0)=x$ and $\dot{\sigma}(t)= T'_{\sigma(t)}$, for large $r(x)$. We can push it into a smooth loop $\tau:=\pi \circ \sigma$ on $\Rl^m$. Observe $L=L'$ is well defined since it is the injectivity radius at infinity. The point is $\tau(L)=\tau(0)$, so that Taylor formula provides   
$
L \dot{\tau}(0) = \int_0^L t \ddot{\tau}(t) dt.
$  
If we decompose $T'$ as the sum of its $h$-horizontal part $H$ and $h$-vertical part $W$, we deduce
\begin{equation}\label{H1}
\abs{H} = \abs{\pi_* T'} = \abs{\dot{\tau}(0)} \leq \int_0^L  \abs{\ddot{\tau}(t)} dt.
\end{equation}
We therefore need an estimate on
\begin{equation}\label{H2}
\abs{\ddot{\tau}} = \abs{D_{\pi_*{T'}} \pi_*{T'}} = \abs{\pi_*(\nabla^h_{H} H)} 
= \abs{\pi_*(\nabla^h_{H} T') - \pi_*(\nabla^h_{H} W)}.
\end{equation}
Since for any $i$, we have 
$$
\abs{h(\nabla^h_{H} W,X_i)} = \abs{h(W,\nabla^h_{H} X_i)} = \frac{\abs{W} \abs{d\eta(H,X_i)}}{2}
\leq c \; r^{1-m},
$$ 
we get $\abs{\ddot{\tau}} \leq c \abs{\nabla^h T'} + c r^{1-m} \leq c \; r^{1-m}$. With (\ref{H1}) and (\ref{H2}), we obtain $H = \OO(r^{1-m})$. Therefore, up to an error term of order $r^{1-m}$, the vector fields $T$ and $T'$ are colinear. Since they have the same $h$-norm up to $\OO(r^{2-m})$, we may assume $T = T' + \OO(r^{2-m})$ (changing $T'$ into $-T'$ if necessary). This ensures we can find convenient orthonormal frames for both metrics: given the orthonormal frame $(X_1,\dots,X_m,T)$ for $h$, we can choose the orientation of the vector field $T'$ so that $T'$ and $T$ are close to each other; if we average the vectors $\pi'_* X_i$ along the fibers of $\pi'$, we get a frame field $(\check{X}_1,\dots,\check{X}_m)$ on the base $\Rl^m$ of $\pi'$ that is asymptotic to an orthonormal frame $(v_1,\dots,v_m)$ of $\Rl^m$ ; this frame induces coordinate functions $y'_1,\dots, y'_m$ on the base and we lift them as functions $x'_i=\pi'^* y'_i$ on $M$. The corresponding frame $(X'_1,\dots,X'_m,T')$ is the one we work with (we can because formula (\ref{formule1}) implies the mass $\Tr \mathcal{Q}_{g,h'}$ does not depend on the choice of $(X'_1,\dots,X'_m,T')$). 

To complete the second step of the proof, we observe the following formula (cf. (\ref{ipp})):
$$
\int_\Omega \left[ \abs{\nabla \alpha'_k}^2 + (\Ric (\alpha'_k) , \alpha'_k) \right] dvol 
=  \int_{\partial \Omega} * \zeta_{\alpha'_k,\alpha'_k},
$$
valid for any domain $\Omega$. The asymptotic of $g$, together with $\Ric_g \in L^1$ and lemma (\ref{constructharm}) ensure the following integral is convergent:
$$
\int_M \left[ \abs{\nabla \alpha'_k}^2 + (\Ric (\alpha'_k) , \alpha'_k) \right] dvol.
$$  
We can therefore express the mass with respect to $h'$ as
$$
\mu_{g,h'} = \frac{1}{\omega_{m} L} \sum_{k=1}^m 
\int_M \left[ \abs{\nabla \alpha'_k}^2 + (\Ric (\alpha'_k) , \alpha'_k) \right] dvol,
$$
which in turn justifies the formula
\begin{eqnarray*}
\mu_{g,h'} = \frac{1}{\omega_{m} L}
\lim_{R \To \infty} \int_{\partial B_R} *_{h'} \; \Big{(}   \left[
 X'_j \cdot g(X'_i,X'_j)  -  X'_i \cdot g(X'_j,X'_j)  \right] dx'_i 
 - \frac{d\left( g(T',T') \right)}{2}  \Big{)},
\end{eqnarray*}
where $\partial B_R = r^{-1}(R)$ is defined with respect to $h$.

Let us turn to the third step, which consists in proving that
\begin{eqnarray*}
& &\lim_{R \To \infty} \int_{\partial B_R} *_{h'} \; \Big{(}   \left[
 X'_j \cdot g(X'_i,X'_j)  -  X'_i \cdot g(X'_j,X'_j)  \right] dx'_i
 - \frac{d\left( g(T',T') \right)}{2}  \Big{)} \\
&=&\lim_{R \To \infty} \int_{\partial B_R} *_{h} \; \Big{(}   \left[
 X_j \cdot g(X_i,X_j)  -  X_i \cdot g(X_j,X_j)  \right] dx_i
 - \frac{d\left( g(T,T) \right)}{2}  \Big{)}.
\end{eqnarray*}
We first write
\begin{eqnarray*}
& & d (g(T',T')) - d (g(T,T)) \\
&=& (\nabla^{h'} g)(T',T') - (\nabla^{h'} g)(T,T) + 2 g(\nabla^{h'} T,T) - 2 g(\nabla^{h'} T',T') \\
&=& (\nabla^{h'} g)(T,T) - (\nabla^{h'} g)(T,T) + 2 h'(\nabla^{h'} T,T') - 2 h'(\nabla^{h'} T',T') + \OO(r^{3-2m}) \\
&=& 2 \; \eta'(\nabla^{h'} T) + \OO(r^{3-2m}),
\end{eqnarray*}
so as to obtain
$$
\lim_{R \To \infty} \int_{\partial B_R} \left[ *_{h'} d (g(T',T')) - *_h d (g(T,T)) \right] 
= 2 \lim_{R \To \infty} \int_{\partial B_R} *_{h} \eta'(\nabla^{h'} T). 
$$
Using $[X_i,T]=0$ and (\ref{connect}), we find
\begin{eqnarray*}
*_{h} \eta'(\nabla^{h'} T) 
&=& \eta'(\nabla^{h'}_{X_i} T) *_{h} dx_i  + \eta'(\nabla^{h'}_{T} T) *_{h} \eta \\
&=& \eta'(\nabla^{h'}_T X_i) *_{h} dx_i + \eta'(\nabla^{h'}_{T} T) *_{h} \eta \\
&=& T \cdot \eta'(X_i) *_{h} dx_i + \eta'(\nabla^{h'}_{T} T) *_{h} \eta  + \OO(r^{3-2m}).
\end{eqnarray*}
The point is:
$$
d\left(*_{h}(\eta' \wedge \eta)\right) 
= T \cdot \eta'(X_i) *_{h} dx_i - X_i \cdot \eta'(X_i) *_{h} \eta.
$$
Since the terms involving $*\eta$ integrate to zero on $B_R$, we can use Stokes formula
to deduce
$$
\lim_{R \To \infty} \int_{\partial B_R} *_{h'} d (g(T',T'))  
= \lim_{R \To \infty} \int_{\partial B_R} *_h d (g(T,T)).  
$$
To tackle the remaining term, we expand it as follows:
\begin{eqnarray*}
& &\left[X'_j \cdot g(X'_i,X'_j)  -  X'_i \cdot g(X'_j,X'_j)  \right] - \left[X_j \cdot g(X_i,X_j)  -  X_i \cdot g(X_j,X_j)  \right]\\
&=& \Big{[}  (\nabla^{h'}_{X'_j} g)(X'_i,X'_j)  +  g(\nabla^{h'}_{X'_j} X'_i,X'_j) + g(X'_i,\nabla^{h'}_{X'_j} X'_j)
- (\nabla^{h'}_{X'_i} g)(X'_j,X'_j) \\
& &- 2g(\nabla^{h'}_{X'_i} X'_j,X'_j) \Big{]} - \Big{[} (\nabla^{h'}_{X_j} g)(X_i,X_j)  +  g(\nabla^{h'}_{X_j} X_i,X_j) \\
& &+ g(X_i,\nabla^{h'}_{X_j} X_j) - (\nabla^{h'}_{X_i} g)(X_j,X_j) - 2g(\nabla^{h'}_{X_i} X_j,X_j) \Big{]}.
\end{eqnarray*}
Using the closeness of the frames and $\nabla^{h'}g = \OO(r^{1-m})$, one can see that the contribution of the terms involving 
$\nabla^{h'}g$ is of order $r^{3-2m}$, hence negligible. Besides, using the closeness of $g$ to $h'$ and (\ref{connect}), we have $g(\nabla^{h'}_{X'_i} X'_j,X'_k) = \OO(r^{3-2m})$, which ensures many terms above are lower order terms. 
Another simplification comes from the fact that the commutator $[X_i,X_j]$ is $\pi$-vertical and of order $r^{1-m}$: it implies 
$\nabla^{h'}_{X_i} X_j$ and $\nabla^{h'}_{X_j} X_i$ only differ by a lower order term. All in all, we find  
\begin{eqnarray*}
& &\left[X'_j \cdot g(X'_i,X'_j)  -  X'_i \cdot g(X'_j,X'_j)  \right] - \left[X_j \cdot g(X_i,X_j)  -  X_i \cdot g(X_j,X_j)  \right]\\
&=& g(\nabla^{h'}_{X_i} X_j,X_j) - g(X_i,\nabla^{h'}_{X_j} X_j) + \OO(r^{3-2m}).
\end{eqnarray*}
We now introduce the $(m-1)$-form $h'(X_j,X'_i) *_h(dx_i \wedge dx_j)$ and compute its differential, up to lower order terms:
\begin{eqnarray*}
  d \left( h'(X_j,X'_i) *_h(dx_i \wedge dx_j) \right)
= g(\nabla^{h'}_{X_j}X_j,X_i) *_h dx_i - g(\nabla^{h'}_{X_i}X_j,X_i)*_h dx_j
+ \OO(r^{3-2m}).
\end{eqnarray*}
After changing $\nabla^{h'}_{X_i}X_j$ into $\nabla^{h'}_{X_i}X_j$ (which costs only a lower order term) and switching summation 
indices $i$ and $j$, only in the second term, we eventually find  
\begin{eqnarray*}
  d \left( h'(X_j,X'_i) *_h(dx_i \wedge dx_j) \right)
= g(\nabla^{h'}_{X_j}X_j,X_i) *_h dx_i - g(\nabla^{h'}_{X_i}X_j,X_j)*_h dx_i
+ \OO(r^{3-2m}).
\end{eqnarray*}
This computation can be combined with Stokes formula to ensure
\begin{eqnarray*}
& &\lim_{R \To \infty} \int_{\partial B_R} *_{h'} \left[X'_j \cdot g(X'_i,X'_j)  -  X'_i \cdot g(X'_j,X'_j)  \right] \\
&=& \lim_{R \To \infty} \int_{\partial B_R} *_h \left[X_j \cdot g(X_i,X_j)  -  X_i \cdot g(X_j,X_j)  \right]. 
\end{eqnarray*}
We have proved: $\mu_{g,h'}^{GB} = \mu_{g,h}^{GB}$.
\endproof

\subsection{The spin case, with a trivial fibration.}

In this paragraph, we would like to explain a spin analogue to the construction above. The setting is a complete \emph{spin} manifold $(M^{m+1},g)$, $m \geq 3$, such that for some compact subset $K$, $M \backslash K$ is diffeomorphic to $\Rl^m \backslash \mathbb{B}^m \times \Sph^1$, which we endow with the standard flat metric $h_0 = dx^2 + dt^2$. We assume the metric $g$ is asymptotic to $h_0$ in the following sense:
$$
g = h_0 + \OO(r^{2-m}), \quad D g = \OO(r^{1-m}) \quad
\text{and} \quad D^2 g = \OO(r^{-m}).
$$
We furthermore assume the spin structure of $M$ coincides outside $K$ with the \emph{trivial} spin structure on $\Rl^m \backslash \mathbb{B}^m \times \Sph^1$.

We work on the spinor bundle $\Sigma M := \Sigma^g M (=:E)$ corresponding to the metric $g$ and we enow it with the pullback of the Levi-Civita connection. Then $\DD$ is the standard Dirac operator $\dirac$ and $\RR$ is the multiplication by $\frac{1}{4} \Scal$ (\cite{LM}). 

Outside $K$, $\Sigma M$ is \emph{not} the bundle corresponding to $h_0$, but we can identify them in a natural way. To make it precise, we denote by $P$ the unique $g$-symmetric section of $\End TM$ such that $h_0= g(P.,P.)$. One can also see $P$ as a natural bijection between the principal bundles of orthonormal frames so that it lifts as an identification of the spin bundles
$$
P :\; Spin(T(M \backslash K),h_0) \; \stackrel{\sim}{\To} \; Spin(T(M \backslash K),g)
$$
and therefore enables to identify $\Sigma^{h_0} (M \backslash K)$ and $\Sigma^g (M \backslash K) = \Sigma (M \backslash K)$. The Levi-Civita connection $D$ of $h_0$ induces a flat metric connection $\nabla^{euc}$ on $(TM,g)$, given by the formula $\nabla^{euc}_X Y = P D_X (P^{-1} Y)$. Since it is a metric connection, it induces a metric connection $\nabla^{euc}$ on the spinor bundle $\Sigma^g M = \Sigma M$. To sum up, we have three connections on $TM$: $\nabla^g$ is metric for $g$ and torsionless; $\nabla^{euc}$ is metric for $g$, is flat but has torsion; $D$ is metric for $h_0$, flat and torsionless. Only two of them lift to $\Sigma M$: $\nabla^g$ and $\nabla^{euc}$.

Now, the spinor bundle $\Sigma (M \backslash K)$ is trivial and $\nabla^{euc}$-flat, so we can find a unit $\nabla^{euc}$-parallel spinor field $\alpha_0$. If $\epsilon$ is a small positive number, we can adapt lemma \ref{constructharm}.

\begin{lem}\label{constructharmdir}
There is a spinor field $\alpha := \widetilde{\alpha} + \beta$ such that $\dirac \alpha =0$, with $\widetilde{\alpha} = \alpha_0$ outside a compact set, $\beta = \OO(r^{2-m+\epsilon})$ and $\nabla \beta \in L^2_{1 - \frac{m}{2} + \epsilon}$.
\end{lem}

\proof
If $\chi$ is a convenient truncature function $\chi$, we can see $\widetilde{\alpha}_0 := \chi \alpha_0$ as a section of $\Sigma M$. Set $\gamma := - \dirac (\chi \alpha_0)$. Since $\gamma = \OO(r^{1-m})$, it belongs to $L^2_{1 -\frac{m}{2}+\epsilon}$. From analysis in weighted spaces (corollary \ref{surjdirac}), we obtain a solution $\sigma \in H^2_{3 -\frac{m}{2}+\epsilon}$ of the equation $\dirac^2 \sigma =\gamma$. Put $\beta = \dirac \sigma$. The estimates on $\beta$ follow as in lemma \ref{constructharm}.  
\endproof

As in the Gauss-Bonnet case, we can use formula \ref{ipp} to find
\begin{equation}\label{formuledir}
\liminf_{R \To \infty} \int_{B_R} \Big{[} \abs{\nabla^g \alpha}^2 + \Scal_g \abs{\alpha}^2_g \Big{]} dvol \leq \limsup_{R \To \infty}\int_{\partial B_R} * \zeta_{\alpha_0,\alpha_0}. 
\end{equation}
and we are interested in computing the right-hand side. 
\begin{lem}
We have 
$
\zeta_{\alpha_0,\alpha_0} = -\frac{1}{4} \left( d \Tr_{h_0} g + \Div_{h_0} g \right) + \OO(r^{3-2m}).
$
\end{lem}

\proof The proof is by now standard. We consider the frame field $(\partial_a)_a := (\partial_1,\dots,\partial_m,\partial_t)$. It is orthonormal for the model metric $h_0$. Putting $e_a := P \partial_a$, we obtain an orthonormal frame field for $g$. We need to understand
$
\zeta_{\alpha_0,\alpha_0} 
= \frac{1}{2} g([e_a \cdot ,e_b \cdot ]  \nabla^g_{e_b} \alpha_0 , \alpha_0 ) e^a.
$
Since $\alpha_0$ is $\nabla^{euc}$-parallel, with $\omega_{c d} := g(\nabla^g e_c, e_d)$ and $\nabla^g - \nabla^{euc} = \frac{1}{8} \omega_{c d} [e_c \cdot , e_d \cdot]$  (\cite{LM}), we obtain
\begin{equation}\label{exprzeta0}
\zeta_{\alpha_0,\alpha_0} 
= \frac{1}{16} \omega_{c d}(e_b) g([e_a \cdot ,e_b \cdot ] [e_c \cdot , e_d \cdot ] \alpha_0 , \alpha_0 ) e^a.
\end{equation}
In order to compute the connection one-form, we resort to Koszul formula:
$$2 \omega_{c d}(e_b) = - g(e_b,[e_c,e_d]) - g(e_c,[e_b,e_d]) + g(e_d,[e_b,e_c]).
$$
Expanding the brackets thanks to the torsionless connection $D$, one finds:
\begin{equation}\label{omega}
2 \omega_{c d}(e_b) 
= -(D_{e_d} g)(e_b,e_c) + (D_{e_c}g)(e_b,e_d) 
+ g(D_{e_b} e_c,e_d)- g(D_{e_b} e_d,e_c).
\end{equation}
Let us denote by $H$ the $g$-symmetric endormorphism such that $g - h_0 =g(H.,.)$. From $P^2 = I - H$ and $e_c = P \partial_c$, we get $D_{e_b} e_c = - \frac{D_{e_b}H}{2} \partial_c + \OO(r^{3-2m})$. And since $H$ is $g$-symmetric, we have
$
g((D_{e_b}H) \partial_c, \partial_d) = g((D_{e_b}H) \partial_d, \partial_c) + \OO(r^{3-2m}).
$
So we deduce $g(D_{e_b} e_c,e_d)- g(D_{e_b} e_d,e_c) = \OO(r^{3-2m})$.
Plugging this into (\ref{omega}) and then (\ref{exprzeta0}), we see that $\zeta_{\alpha_0,\alpha_0}$ can be approximated by
\begin{eqnarray*}
& & \frac{1}{32} \left( \partial_c g_{b d} - \partial_d g_{b c} \right) 
g([e_a \cdot ,e_b \cdot ] [e_c \cdot, e_d \cdot] \alpha_0 , \alpha_0 ) e^a + \OO(r^{3-2m}) \\
&=& \frac{1}{16}  \partial_c g_{b d} \; g([e_a \cdot ,e_b \cdot ] [e_c \cdot, e_d \cdot] \alpha_0 , \alpha_0 ) e^a + \OO(r^{3-2m}) \\
&=& \frac{1}{4}  \partial_c g_{b d} \; g(\left( \delta_{a b} \delta_{c d} + \delta_{a b} e_c e_d + \delta_{c d} e_a e_b + e_a e_b e_c e_d \right) \cdot \alpha_0 , \alpha_0 ) e^a + \OO(r^{3-2m}).
\end{eqnarray*}
In view of (\ref{comprod}) and (\ref{comzero}), every term like
$g(e_a \cdot e_b \cdot \alpha_0 , \alpha_0 )$ can be replaced by $-\delta_{a b} \abs{\alpha_0}^2$. Moreover, using (\ref{comprod}), we can write
$$
\partial_c g_{b d} \; g(e_a e_b e_c e_d \cdot \alpha_0 , \alpha_0 )
= - \frac{1}{2} \partial_c g_{b d} \; g(e_a [e_b, e_d] e_c \cdot \alpha_0 , \alpha_0 )
+ \partial_c g_{b d} \; g(\delta_{b d} e_a  e_c \cdot \alpha_0 , \alpha_0 ).
$$
Since $\partial_c g_{b d}$ is symmetric with respect to $b$ and $d$ whereas $[e_b,e_d]$ is antisymmetric, the first term vanishes. These observations lead to :
\begin{eqnarray*}
\zeta_{\alpha_0,\alpha_0} &=& \frac{1}{4}  \partial_c g_{b d} \; \left( \delta_{a b} \delta_{c d} - \delta_{a b} \delta_{c d}  - \delta_{c d} \delta_{a b} -  \delta_{b d} \delta_{a c} \right) \abs{\alpha_0}^2 e^a + \OO(r^{3-2m}) \\
&=& \frac{1}{4} \; \left( - \partial_c g_{a c} - \partial_a g_{b b} \right) e^a + \OO(r^{3-2m}),
\end{eqnarray*}
hence the lemma. 
\endproof 

The corresponding positive mass theorem involves
\begin{eqnarray}
\mu_{g,h}^{D} = -\frac{1}{\omega_{m} L}
\limsup_{R \To \infty} \int_{\partial B_R} *_{h} \; \Big{(} \Div_{h_0} g + d \Tr_{h_0} g \Big{)}.
\end{eqnarray}
\begin{thm}
Let $(M^{m+1},g)$, $m \geq 3$, be a complete \emph{spin} manifold with nonnegative scalar curvature. We assume there is a compact set $K$ and a spin preserving diffeomorphism between $M \backslash K$ and $\Rl^{m} \backslash \Ball^m \times \Sph^1$ such that  
$$
g = g_{\Rl^m \times \Sph^1} + \OO(r^{2-m}), \quad D g = \OO(r^{1-m}) \quad
\text{and} \quad D^{2} g = \OO(r^{-m}).
$$
Then $\mu_{g,h}^D$ is nonnegative and vanishes exactly when $(M,g)$ is isometric to the standard
$\Rl^m \times \Sph^1$.
\end{thm}

\proof
Since $\Scal_g \geq 0$, formula \ref{formuledir} leads to  
$
\int_{M} \abs{\nabla^g \alpha}^2 dvol 
\leq  \frac{\omega_{m} L}{4} \mu_{g,h}^D, 
$
hence the nonnegativity of $\mu_{g,h}^D$. When $\mu_{g,h}^D$ vanishes, every constant spinor $\alpha_0$ in the model gives rise to a harmonic and parallel spinor field $\alpha$ that is asymptotic to $\alpha_0$. This makes it possible to produce a parallel trivialization of the spinor bundle. It follows that $(M,g)$ has trivial holonomy (\cite{MS}) and we can conclude as in theorem \ref{GBmass}.
\endproof

Finally, the proof of proposition \ref{welldef} implies
\begin{prop}
In this setting, provided the scalar curvature is in $L^1$, $\mu_{g,h}^D$ does not depend on $h$ but only on $g$.
\end{prop}

\vskip 0.5cm

\section{Examples.}

\subsection{Schwarzschild metrics.}

These are complete Ricci flat metrics on $\Rl^2 \times \Sph^{n-2}$, $n\geq 4$,
given by the formula:
$$
g_{\gamma} = d\rho^2 
+ F_\gamma(\rho)^2 d\theta^2 + G_\gamma(\rho)^2 d\omega^2 .
$$   
$\rho$, $\theta$ are polar coordinates on the $\Rl^2$ factor, $d\omega^2$ is the standard metric on $\Sph^{n-2}$, $F_\gamma$ and $G_\gamma$ are smooth functions defined by
$$
G_\gamma'(\rho) = \sqrt{1 - \left(\frac{\gamma}{G_\gamma}\right)^{n-3}}, \quad
G(0)=\gamma
\quad \text{and} \quad  
F_\gamma = \frac{2 \gamma}{n-3} \sqrt{1 - \left(\frac{\gamma}{G_\gamma}\right)^{n-3}},
$$
for some positive parameter $\gamma$. $G_\gamma$ increases from $\gamma$ to $\infty$ and $G_\gamma \sim \rho$ at infinity ; $F_\gamma$ increases from $0$ to $\frac{2 \gamma}{n-3}$ and $F_\gamma \sim \rho$ near $0$. Setting $r := G_\gamma(\rho)$ and $t = \frac{2\gamma}{n-3} \theta$, we can write
$$
g_\gamma = \frac{dr^2}{1 - \left( \frac{\gamma}{r} \right)^{n-3}} + r^2 d\omega^2 +  \left[ 1 - \left( \frac{\gamma}{r} \right)^{n-3} \right] dt^2 
$$
In this way, it is apparent that $g_\gamma$ is asymptotic to the flat metric on $\Rl^{n-1} \times \Sph^1$, with circle length equal to  $L := \frac{4 \pi \gamma}{n-3}$ at infinity. Observe $M = \Rl^2 \times \Sph^{n-2}$ is spin (with a unique spin structure, since it is simply connected). 

A fancy way to compute the masses this consists in introducing ``isotropic'' coordinates: we define a new radial coordinate $u$ by
$$
r = u \left[ 1 + \frac{1}{4} \left( \frac{\gamma}{u} \right)^{n-3} \right]^{\frac{2}{n-3}},
$$
so that Schwarzschild metrics are described by 
$$
g_\gamma = \left[ 1 + \frac{1}{4} \left( \frac{\gamma}{u} \right)^{n-3} \right]^{\frac{4}{n-3}} (du^2 + u^2 d\omega^2) + \left[  \frac{1 - \frac{1}{4} \left( \frac{\gamma}{u} \right)^{n-3}}{1 + \frac{1}{4} \left( \frac{\gamma}{u} \right)^{n-3}}   \right]^2 dt^2
$$
for $u> 0$. We can then choose Cartesian coordinates $x_1, \dots, x_{n-1}$ corresponding to the polar coordinates $(u, \omega)$ and keep in mind the first order terms:  
$$
g_\gamma \approx \left[ 1 + \frac{\gamma^{n-3}}{n-3} u^{3-n} \right] dx^2 
+ \left[ 1 - \gamma^{n-3} u^{3-n} \right] dt^2.
$$
With $h_0=dx^2 + dt^2$, this readily provides 
$d\Tr_{h_0} g_\gamma = -2  \gamma^{n-3} u^{2-n} du + \OO(u^{5-2n})$,
$\Div_{h_0} g_\gamma = \gamma^{n-3} u^{2-n} du + \OO(u^{5-2n})$ and 
$d g_\gamma(\partial_t,\partial_t) = (n-3)  \gamma^{n-3} u^{2-n} du + \OO(u^{5-2n})$.
In view of the definitions, we deduce
$$
\mu_{g_\gamma}^{D} = \gamma^{n-3} \quad \text{and} \quad \mu_{g_\gamma}^{GB} = \frac{n-1}{2} \gamma^{n-3} = \frac{n-1}{2} \mu_{g_\gamma}^{D}.
$$
For instance, in dimension $n=4$, the masses reduce to the parameter $\gamma$ (up to a constant), which is basically what we hoped (cf. introduction). 

The mass is therefore positive, but we cannot deduce it from the ``spin'' version of the positive mass theorem: the spin structure at infinity is not the trivial spin structure! Indeed, the spin structure on the asymptotic $\Sph^1$ comes from the spin structure of the unit circle in the $\Rl^2$ factor of $\Rl^2 \times \Sph^{n-2}$. This spin structure is therefore induced by the (unique and trivial) spin structure on the unit disk: it is the non trivial spin structure on the circle (see \cite{Bar}). In the next paragraph, we will give examples showing that this non trivial spin structure at infinity really allows negative mass in nonnegative scalar curvature. But here, since the Ricci curvature of Schwarzschild metrics is nonnegative, we can use the Gauss-Bonnet point of view and explain the positivity of the mass: it is a consequence of $\Ric \geq 0$ and not of $\Scal \geq 0$. 

\subsection{Reissner-Nordstr\"om metrics.}

In dimension $4$, one can include the Schwarzschild metric $g_\gamma$ in a broader family of complete scalar flat metrics on $\Rl^2 \times \Sph^2$ that are asymptotic to $\Rl^3 \times \Sph^1$, with the non trivial spin structure. These Reissner-Nordstr\"om metrics \cite{Dai,BrH,CJ,PK} are given by the same ansatz $g = d\rho^2 + F(\rho)^2 d\theta^2 + G(\rho)^2 d\omega^2$, with
$$
G' = \sqrt{1 - \frac{2m}{G} - \frac{q^2}{G^2} }
\quad \text{and} \quad  
F = \frac{G_0^2}{G_0 - m} 
\sqrt{1 - \frac{2m}{G} - \frac{q^2}{G^2} }
$$
and $G(0) := G_0 = m + \sqrt{m^2 + q^2}$. The behaviour of $G$ and $F$ is similar to the Schwarzschild analogue. The new feature is $m$ can be chosen negative: the metric is then still complete, provided $q$ is nonzero. Setting $r := G(\rho)$ and $t = \frac{G_0^2}{G_0 - m}  \theta$, we can obtain a more familiar formula:
$$
g = \frac{dr^2}{1 - \frac{2m}{r} - \frac{q^2}{r^2}} +  \left[ 1 -  \frac{2m}{r} - \frac{q^2}{r^2}  \right] dt^2 + r^2 d\omega^2.
$$
The formulas for doubly-warped products \cite{Pet} (or geometric arguments as in \cite{Bes}, 3.F) make it possible to compute the curvature. The eigenvalues of the Ricci tensor are $\frac{q^2}{G^4}$, along $\partial_\rho$ and $\partial_t$, and $-\frac{q^2}{G^4}$, along the $\Sph^2$ factor. It therefore has no sign but the scalar curvature vanishes. To compute the mass
$\mu_g^{D}$, we can again use isotropic coordinates. The new radial coordinate 
$u$ is given by
$
r = u \left[ 1 + \frac{m}{u} + \frac{m^2 + q^2}{4u^2}  \right] 
$
and we find
$$
g = \left[ 1 +  \frac{m}{u} + \frac{m^2 + q^2}{4u^2}  \right]^2 (du^2 + u^2 d\omega^2) 
+ \left[  \frac{1 - \frac{m^2 + q^2}{4u^2} }{1 + \frac{m}{u} + \frac{m^2 + q^2}{4u^2} }   \right]^2 dt^2.
$$ 
Comparing to the Schwarzshild formula, one can see that the asymptotic is the same up to  $\OO(u^{-2})$, so the mass $\mu_g^{D}$ is $2m$. Since we can choose $m$ negative, this yields a whole family of metrics with negative mass!

\subsection{Taub-NUT metrics.}

The Taub-NUT metric is the basic non trivial example of ALF gravitational instanton. In particular, it is Ricci flat. For a precise description, 
the reader is referred to \cite{Leb}. Basically, it is a complete metric on $\Rl^4$ 
which is adapted at infinity to the Hopf fibration: it can be written outside one point as
$$
g = \left( 1 + \frac{2m}{r} \right) dx^2  + \frac{1}{1 + \frac{2m}{r}} \eta_1^2.
$$
for some positive parameter $m$. Here, $\eta_1$ is the standard contact form on the three-spheres
(cf. section 1). 
With $h_1 = dx^2 + \eta_1^2$, we can compute:
$ d \Tr_{h_1} g = -\frac{4m dr}{r^2}  + \OO(r^{-3})$, 
$ d g(T,T) = \frac{2m dr}{r^2}  + \OO(r^{-3})$. With (\ref{connect}), we also find 
$\Div_{h_1} g = \frac{2m dr}{r^2}  + \OO(r^{-3})$.
This leads to $\mu_{g}^{GB}= 3m$ (which is similar to the Schwarzschild case).

The same computations apply to the multi-Taub-NUT metrics (\cite{Leb}). These metrics 
are again hyperk\"ahler hence Ricci flat and their asymptotic is that of $h_k$. Their 
mass $\mu^{GB}$ is easily seen to be equal to $3km$.



\end{document}

%% file: expcrit.pstex_t
\begin{picture}(0,0)%
\includegraphics{expcrit.pstex}%
\end{picture}%
\setlength{\unitlength}{3947sp}%
\begingroup\makeatletter\ifx\SetFigFont\undefined%
\gdef\SetFigFont#1#2#3#4#5{%
  \reset@font\fontsize{#1}{#2pt}%
  \fontfamily{#3}\fontseries{#4}\fontshape{#5}%
  \selectfont}%
\fi\endgroup%
\begin{picture}(6012,664)(976,-719)
\put(2776,-211){\makebox(0,0)[lb]{\smash{{\SetFigFont{12}{14.4}{\familydefault}{\mddefault}{\updefault}{\color[rgb]{0,0,0}$2-\frac{m}{2}$}%
}}}}
\put(2176,-661){\makebox(0,0)[lb]{\smash{{\SetFigFont{12}{14.4}{\familydefault}{\mddefault}{\updefault}{\color[rgb]{0,0,0}$1-\frac{m}{2}$}%
}}}}
\put(5251,-661){\makebox(0,0)[lb]{\smash{{\SetFigFont{12}{14.4}{\familydefault}{\mddefault}{\updefault}{\color[rgb]{0,0,0}$\frac{m}{2}+1$}%
}}}}
\put(4576,-211){\makebox(0,0)[lb]{\smash{{\SetFigFont{12}{14.4}{\familydefault}{\mddefault}{\updefault}{\color[rgb]{0,0,0}$\frac{m}{2}$}%
}}}}
\put(976,-661){\makebox(0,0)[lb]{\smash{{\SetFigFont{12}{14.4}{\familydefault}{\mddefault}{\updefault}{\color[rgb]{0,0,0}$\dots$}%
}}}}
\put(6601,-661){\makebox(0,0)[lb]{\smash{{\SetFigFont{12}{14.4}{\familydefault}{\mddefault}{\updefault}{\color[rgb]{0,0,0}$\dots$}%
}}}}
\put(5851,-211){\makebox(0,0)[lb]{\smash{{\SetFigFont{12}{14.4}{\familydefault}{\mddefault}{\updefault}{\color[rgb]{0,0,0}$\frac{m}{2}+2$}%
}}}}
\put(1576,-211){\makebox(0,0)[lb]{\smash{{\SetFigFont{12}{14.4}{\familydefault}{\mddefault}{\updefault}{\color[rgb]{0,0,0}$-\frac{m}{2}$}%
}}}}
\end{picture}%